\documentclass[12pt, reqno]{amsart}

\usepackage{philstyle}

\def\cF{{\mathcal{F}}}

\def\tD{\tilde{D}}

\def\R{{\mathbb R}}
\def\T{{\mathbb T}}
\def\Z{{\mathbb Z}}
\def\N{{\mathbb N}}

\def\C{{\mathbb C}}

\def\raw{\rightarrow}

\def\I{^{-1}}

\def\td{\tilde{d}}

\def\tx{{\tilde{x}}}
\def\ty{{\tilde{y}}}
\def\tY{{\tilde{Y}}}

\def\tX{{\tilde{X}}}

\def\tz{{\tilde{z}}}

\def\tf{\tilde{f}}
\def\tfk{\tilde{f}^{(k)}}
\def\tF{\tilde{F}}
\def\tg{\tilde{g}}

\def\tphi{\tilde{\phi}}

\def\th{{\tilde{h}}}

\def\tiota{{\tilde{\iota}}}

\def\tMA{\tilde{M}_a}
\def\tXA{\tilde{X}_a}
\def\tMu{\tilde{M}_u}
\def\tXu{\tilde{X}_u}
\def\tM{\tilde{M}}
\def\tMa{\tilde{M}_a}

\def\vn{{\vec{n}}}
\def\vm{{\vec{m}}}
\def\deltan{\delta_\vn}
\def\tomega{\tilde{\omega}}
\def\talpha{\tilde{\alpha}}
\def\tbeta{\tilde{\beta}}
\def\tsigma{\tilde{\sigma}}
\def\tgamma{\tilde{\gamma}}

\def\htop{h_{top}}
\def\homeo{homeomorphism}

\DeclareMathOperator{\diam}{diam}

\DeclareMathOperator{\NC}{NC}

\DeclareMathOperator{\BF}{BF}

\DeclareMathOperator{\Fix}{Fix}
\DeclareMathOperator{\fix}{Fix}
\DeclareMathOperator{\trace}{trace}

\def\directlim{{\varinjlim}}

\def\hsigma{{\hat{\sigma}}}

\def\vzero{\vec{0}}

\def\vv{{\vec{v}}}
\def\vw{{\vec{w}}}

\def\floor#1{\lfloor #1 \rfloor}
\def\mc#1{\textit{#1}}
\def\mycol#1{\textit{#1}}
\def\pA{pseudoAnosov}

\def\Inv{^{-1}}

\def\titf{\tiota\tf}
\def\titfk{\tiota\tfk}
\def\titx{\tiota(\tx)}

\usepackage{amssymb,amsmath,amscd,xspace}
\usepackage{enumerate}
\usepackage{graphicx}

\graphicspath{ {./figures/} }

\begin{document}

\title{Dynamical displacements, persistence and semiconjugacies}

\author{Philip Boyland}
\address{Department of
    Mathematics\\University of Florida\\372 Little Hall\\Gainesville\\
    FL 32611-8105, USA}
\email{boyland@ufl.edu}

\begin{abstract}
This survey gives a unified treatment of topics from Abelian and non-Abelian
Nielsen Theory integrated with the  semiconjugacy theorems of Franks and Handel.
The main focus is to develop an analog of the rotation set that is 
valid when the dynamics are not isotopic to the identity and to connect
this theory to the dynamical persistence under homotopy/isotopy
intrinsic in the theorems of Franks and Handel. For this
dynamical persistence,  expansion/hyperbolicity at some scale is
 essential.
\end{abstract}

\maketitle

\section{Introduction}

When confronted with complicated systems, a standard mathematical
strategy is to construct simpler invariants which reveal at least
part of the structure. In dynamical systems on manifolds one implementation
of this strategy begins with the rotation vector of an orbit. The direction
of the rotation vector indicates the direction of the asymptotic motion of the
orbit in first homology while  its length indicates the speed. The
rotation vector is computed by lifting the dynamics to the universal
Abelian cover  and taking a Birkhoff average of the displacements
of the lifted orbit. The rotation vectors of all the orbits
are collected in the desired invariant, the rotation set.

One shortcoming of this strategy is that it only works when 
the dynamics is homotopic to the identity or, more generally, acts like 
the identity on first homology. So if the dynamics is generated
by $f:M\raw M$, the method requires $A(f) := (f_1)_*:H_1(M)\raw H_1(M)$
be $A(f) = I$. The reasons for this restriction are at least
 twofold. First,
when $A(f)\not= I$ the displacement of an orbit in the universal
Abelian cover $\tMA$ depends
on the choice of the lift. Second, if $A(f)$ has eigenvalues outside
the unit circle, orbits will be running off to infinity at an
exponential rate and a simple Birkhoff average won't work.
 Thus when $A(f)\not= I$ different methods are needed
and in this paper we present three, which turn out to be equivalent.

For simplicity of exposition let us assume that $M$ is a smooth,
connected orientable manifold with  $H_1(M;\Z) = \Z^b$, torsion free.
The first method starts with periodic points. Defining the
displacement of a lifted $k$-periodic point
 in the $k^{th}$-order Bowen-Franks group
$$\BF_k(A(f)) = \frac{\Z^b}{(A(f)^k - I)\Z^b}$$ makes it
independent of the choice of the lift. We shall see in
Section~\ref{toral} that $\BF_k(A)$
can be naturally identified with $\fix(\Phi_A^k)$ where
$\Phi_A:\T^b\raw\T^b$ is the toral endomorphism induced by the
matrix $A(f)$. Passing to a direct limit yields $\BF_\infty(A(f))$
which corresponds to the collection of all periodic points
of $\Phi_A$. Thus the displacement of a periodic point
associates it with a unique periodic point of $\Phi_A$.

To work with general, non-periodic orbits we require more conditions
on the action $A(f)$ of $f$  on first homology. Let us assume
it is $H_1$-expanding which means all its eigenvalues are outside
the unit circle ($H_1$-hyperbolic is similar).
The Abel-Jacobi-Albanese map allows us to compare lifted orbits of
$f$ to the universal Abelian cover with lifted orbits of $\Phi_A$
in $\R^b$. The $H_1$-expanding hypothesis yields that  lifts of
periodic points of $f$ stay a uniformly bounded distance from the 
lift of the point in $\BF_\infty(A(f))$ to which it corresponds.
This is expressed by saying the two orbits \textit{globally shadow}.
This is a notion that works with general orbits and a variant
of the proof of the shadowing lemma says that each orbit $o(x, f)$
globally shadows an orbit $(y,\Phi_A)$. The assignment $\alpha:x\mapsto y$
yields a semiconjugacy of $f$ acting on $M$ to $\Phi_A$ acting
on the invariant subset $\alpha(M)\subset\T^b$. 

The second approach works directly with the lifted iterates of 
$f$ in $\tMA$ by normalizing $\tf^k$ by the rate $A^{-k}$. The
limit converges and again yields a semiconjugacy $\beta:M\raw\T^b$
of $(M,f)$ into $(\T^b,\phi_A)$. It is shown that $\beta$ is
the same as $\alpha$ obtained from global shadowing.
This method is closest to Franks
original construction of this semiconjugacy in \cite{franks}.

In the third approach a finite segment of a general orbit of $f$
is approximated by the closest periodic orbit of $\Phi_A$. The 
$f$-orbit is assigned  to the limit of these $\Phi_A$-orbits in $\T^b$. The
resulting map turns out to be the same as $\alpha = \beta$ of the
previous paragraphs.
  
The three approaches all yield the
same semiconjugacy
which compares and connects the dynamics of $f$ with those of the
toral endomorphism $\Phi_A$. The image in $\T^b$ captures
all the asymptotic dynamics of $f$ as measured in first homology.

In the special case that $M = \T^b$ 
the semiconjugacy turns out to be onto yielding
the following interpretation. If $f:\T^b\raw \T^b$
is such that $A(f)$ is expanding, then $f$ is semiconjugate
to $\Phi_A$. This means that all the dynamics of $\Phi_A$
are present in any homotopic map. This, and its $H_1$-hyperbolic
version, is the context
in which Franks Theorem is usually seen.

The second main part of the paper studies what can be viewed as a 
non-Abelian analog of Franks Theorem. In this case we work
with the universal cover and so start by defining displacements
as twisted conjugacy classes in $\pi_1(M)$. This is the
realm of classical Nielsen Theory. The semi-conjugacy
theorem in this context, due to Handel \cite{handel}, requires 
metric expansion/hyperbolicity in the universal cover along with conditions
on the periodic orbits of $f$. The prototypical example
of such maps are \pA\ homeomorphisms on compact surfaces.

Most of the mathematics covered here is well established. The novel
features are the integrated approach and perhaps Section~\ref{BFinfinity}
including the direct limit of the Bowen-Franks groups. The
essential sources include Jiang's excellent expositions
of Nielsen Theory \cite{jiang1, jiang2} and the
 papers of  Franks \cite{franks},
Handel \cite{handel}, Shub \cite{shub}, Fathi \cite{fathi},
 and Fried \cite{fried1, fried2}. Band first used the connection of 
shadowing to the Franks map in \cite{band}. Bouwman and Kwapisz
added  Abelian Nielsen equivalence to this connection in  their work on
Hirsch's problem (Section ~\ref{hirsch})\cite{BK}
and the connection was explored in a general context by Gromov \cite{gromov}. 

To simplify the exposition we have
mainly dealt with the case of expansion and just remark on
hyperbolicity. Also, rather than break the flow with
frequent references, we include them with comments in  a
Context subsection. The bibliography here is not exhaustive and
we encourage the interested reader to consult the bibliographies
of the references. Finally, we will use maps on the wedge of two 
circles as examples. While the space is obviously not a smooth manifold
all the methods here work with very minor modification.

\noindent\textbf{Acknowledgements:}
Thanks to Gavin Band, Luca de Cerbo, Alex Dranisnikov and Mark Pollicott
 for useful conversations.
This  work is partially supported by the Simons Foundation 
Award No.663281 granted to the Institute of Mathematics of the Polish 
Academy of Sciences for the years 2021-2023.

\section{The universal Abelian covering space}\label{tMA}
As noted in the introduction one of the goals of this paper to measure the
asymptotic motion of orbits in homology. Working in the manifold
itself often leads to ambiguities  and so 
one unwraps the manifold into a covering space,  
lifts the dynamics to the cover, and 
measures the displacement of orbits there. Since we measuring the
motion of orbits in homology, the appropriate covering space is
 the \textit{universal
Abelian cover}, $\pi:\tMA\raw M$.
This space has no natural algebraic structure
outside of the deck group, so we use a classical technique
to construct a map $\tiota:\tMA\raw \R^b$  and use the vector space structure
of the images to measure orbit progress. 
\subsection{The covering space}
For simplicity of exposition, we let $M$ be a smooth, 
connected, orientable, compact manifold 
with $H_1(M;\Z) \cong \Z^b $ torsion-free with $b>0$. A compact, orientable
surface with negative Euler characteristic 
is a good example to keep in mind. The dynamics will be given by
$f:M\raw M$ continuous, onto and orientation preserving.

The \textit{universal Abelian cover} is built by moding out the universal
cover by the action of the commutator subgroup of the fundamental
group. The deck group of $\tMA$ is $\Z^b $  and the deck transformation
corresponding to $\vn\in\Z^b$ is denoted $\deltan$. A fundamental
property of $\tMA$ is that a closed loop $\gamma:[0,1]\raw M$ with $\gamma(0)
= \gamma(1)$ lifts to a closed loop in $\tMA$ if and only if it
is zero in homology. 

The action
of $f$ on the fundamental group preserves the commutator subgroup, and so every
$f$ lifts to $\tf:\tMA\raw \tMA$. The collection of all possible
lifts of $f$ is given by $\deltan\; \tf$ for all $\vn\in\Z^b$. The
fundamental relationship between lifts and deck transformations is
\begin{equation}\label{fundamental}
\tf\circ \deltan = \delta_{A \vn} \circ\tf
\end{equation}
where  $A = (f_*)_1:H_1(M;\Z)\raw H_1(M;\Z)$ is the induced map on homology.

\begin{notation}
Henceforth the action of $f$ on first homology will be denoted $f_*$
and the linear transformation $A = A(f)$ will always represent that action.
We will not be considering higher degree homology in this paper.
\end{notation}

\subsection{Mapping $\tMA$ into $\R^b$}\label{abeljacobi}
Pick a basis $\Gamma_1, \dots, \Gamma_b$ for $H_1(M;\Z)$ and
a basis of harmonic one-forms $\omega, \dots, \omega_b$ for
$H^1_{DR}(M; \R)$ which
are dual in the sense that
\begin{equation}\label{equi}
\omega_i(\Gamma_j) = \delta_{ij}
\end{equation}
(Kronecker delta). Lift the $\omega_i$ to $\tMA$ as $\tomega_i$
and fix a point $\tx_0\in\tMA$. For $\tx\in\tMA$
choose a  smooth path $\gamma$ from $\tx_0$ to $\tx$ and for
$i = 1, \dots, b$ define $\tiota:\tMA\raw\R$ by
\begin{equation}
\tiota_i(\tx) = \int_{\gamma} \tomega_i.
\end{equation}
Since the forms are closed, this is independent of the choice
of path $\gamma$ and  using
\eqref{equi}, 
\begin{equation}\label{deck}
\tiota\circ \deltan = \tiota + \vn
\end{equation} Thus
$\tiota$ descends to map $\iota:M\raw \T^b$.

The map $\iota$ is often not injective. A trivial
example is when the first betti number $b$ is
less than the dimension of the manifold. In many cases of interest,
for example closed surfaces with negative Euler characteristic, $\tiota$
is injective and is thus an embedding. In these cases, one
can treat $\tMA$ as being inside $\R^b$ and eliminate all
the $\tiota$'s from the succeeding formulas. The formulas
are then much neater and intuitive.

However, using \eqref{deck}, $\tiota$ is always bijective on any
set $\{\deltan \tx : \vn\in \Z^b\}$ for $\tx\in\tMa$. In addition,
$\iota$ always induces a surjective map on first homology: given
$\vn\in\Z^b$ pick any path $\gamma$ in $\tMa$  from $\tx_0$
to $\deltan \tx_0$. Then $\pi\circ \gamma$ is closed curve 
 $\Gamma\subset M$ with $(\iota_*)_1 ([\Gamma]) = \vn$, again
using \eqref{deck}. 

\subsection{formula for $\tf$} 
Define $\tsigma:\tMA\raw \R^b$ by $\tsigma(\tx) = \tiota\tf(\tx) - A\tiota(\tx)$
and so $\tsigma$ is continuous.
Using \eqref{deck} and \eqref{fundamental}, 
$\tsigma(\tx + \vn) = \tsigma(\tx)$
and so $\tsigma$ descends to a function on the compact
space $M$. As a consequence we will
usually write just $\sigma(x)$ for $\tsigma(\pi(\tx))$.
Since $M$ is compact, $\sigma$ is bounded.
 The desired formula for $\titf$ is
then
\begin{equation}\label{desired}
\titf(\tx) = A\tiota(\tx) + \sigma(x).
\end{equation}

\subsection{context} The classical version of the map $\iota$ is
the Abel-Jacobi map on a compact Riemann surface using holomorphic
one-forms. The generalization to higher dimensional Riemannian manifolds
is also sometimes called the generalized
Abel-Jacobi map or else the Albanese
map or even the Abel-Jacobi-Albanese map (\cite{gromov}). Its use was introduced into the study of rotation sets by
Mather \cite{mather}.

\section{The translation and the case  $f_*=I$}\label{translation}
\subsection{the translation}
We use the Abelian cover to measure the motion of a point $\tx$
in homology under a single iterate. 
The \textit{translation} in $\tMA$ that a point $\tx$ experiences under
a lift $\tf$  is
\begin{equation}\label{basic}
\Delta(\tx, \tf) := \titf(\tx) - \tiota(\tx) = (A-I)\tiota(\tx) + \sigma(x)
\end{equation}
using \eqref{desired}. Note that $\Delta$ depends on the choice of
lift of $\tf$ is a simple way: changing $\tf$ to $\deltan\tf$ simply
changes $\Delta$ to $\Delta + \vn$. 
On the other hand, 
if we change the lift of $x$ to $\deltan\tx$
the translation becomes
\begin{equation}\label{depends}
\Delta(\deltan \tx) = \Delta(\tx)  + (A-I)\vn.
\end{equation}
 
The goal here is to define an asymptotic average translation
for a given point $\tx\in\tMa$. We want this to depend just
on the point $x$. The formula
\eqref{depends} indicates exactly how to compensate for
changing a lift of $x$ and we will pursue this in
the next section after discussing an important special case.

\subsection{the case $f_* = I$}
Assume $f$ acts trivially on first homology,
 from \eqref{depends}
$\Delta(\deltan \tx) = \Delta(\tx)$ and further, from \eqref{basic}  
$\Delta(\tx, \tf) = \sigma(x)$ and so we may write the displacement as 
a function of $x\in M$ as   
$\Delta(x) = \sigma(x)$
suppressing the choice of lift of $f$.
Iterating the translation after $n$ iterates yields
$\tiota\tf^n(\tx) - \tiota(\tx) = \sigma(x) + \dots \sigma(f^{n-1}(x))$
with the $\sigma$-terms uniformly bounded. Thus 
$\titf^n(\tx) - \tiota(\tx)$ grows
at most linearly and so it is natural to take the Birkhoff average
and define the rotation vector
\begin{equation}\label{birk}
\rho(x) = \lim_{n\raw\infty}\frac{1}{n} \sum_{i=0}^{n-1} \Delta(f^i(x))
= \lim_{n\raw\infty}\frac{\titf^n(\tx) - \tiota(\tx)}{n}.
\end{equation}
This limit  exists at  almost every point of an invariant
measure by Birkhoff's Pointwise Ergodic Theorem. The rotation set
of $f$, $\rho(f)$, is then the union of all the rotation vectors.

The most common situation studied where  $f_* = I$ is when
$f$ is homotopic or isotopic to the identity, but it is worth
noting that there are interesting examples in which $f$ acts
trivially on first homology but the action on the fundamental group
is quite complex. A good example is given by psuedoAnosov maps
in the Torelli subgroup of the mapping class group of
a compact surface \cite{primer}.  The next example shows this phenomenon
in a map on the wedge of two circles.

\begin{figure}[htbp]
\begin{center}
\includegraphics[width=0.85\textwidth]{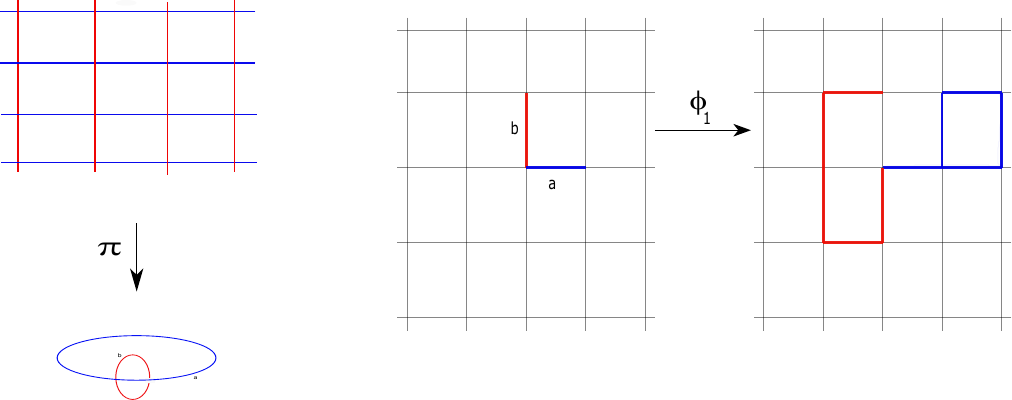}
\end{center}
\caption{Left: The wedge of two circles and its universal Abelian cover.
Right: A fundamental domain and its image under the map $\tphi_1$.}
\label{figure1}
\end{figure}
\subsection{example}\label{eg1}
 Let $X$ be the wedge of two circles as pictured in Figure~\ref{figure1}. The
universal Abelian cover $\tXA$ is $\Z\times\R \cup \R\times\Z$ which
sits naturally in $\R^2$   and we write 
$\tx + \vn$ for $\deltan\tx$.
Consider the homomorphism $\psi_1$ of the free
group on two symbols $F_2 \cong \pi_1(X)$ induced 
by $a\mapsto a a b a\I b\I, b\mapsto b\I a\I b b a$.
  Let $\phi_1:X\raw X$ be the
map that realizes this action in the tightest fashion, so $\phi_1' = 5$
except at the vertex. Lift $\phi_1$ to $\tphi_1:\tX\raw \tX$
which fixes the origin. Figure~\ref{figure1}(right)
 shows the image under $\tphi_1$ of a
 fundamental domain of $\tX$.

\begin{figure}[htbp]
\begin{center}
\includegraphics[width=0.4\textwidth]{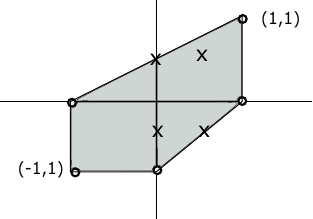}
\end{center}
\caption{The rotation set of the map $\phi_1$.}
\label{rotset}
\end{figure}

While $\phi_1$ is not homotopic to the identity,
 it does act like the identity
on $H_1$ and  so we may compute
its rotation set. Now
 $\tphi_1(\tx+\vn) = \tphi_1(\tx) + \vn$, and so it
suffices to understand the displacements of a fundamental domain.
Letting $t$ be the horizontal deck transformation (adding  $(1,0)$)
 and $r$ the vertical
one (adding $(0,1)$. Using Figure~\ref{figure1}
 we compute the various displacements of
a lift of each of the two circles. This information is encoded
in the transition matrix
\begin{equation*}
M = \begin{pmatrix} 1 + t + tr&  t + t^2 \\ r\Inv t\Inv + t\I  r &
r\I + t\I + r\Inv t\Inv.
\end{pmatrix}  
\end{equation*}
where, for example, $M(1,1)$ says that there are transitions $a\raw a$
which translate by $(0,0), (1,0)$, and $(1,1)$ in the cover.
Next observe that any recurrent allowable sequence of transitions 
can be constructed by
taking the concatenation of minimal loops, i.e.\ cycles of allowable
transitions which make each transition at most once.
Concatenating loops corresponds
to taking the Farey sum of their rotation vectors. As a consequence,
the rotation set is the convex hull of the minimal loops. In this
example the minimal loops are all length one or two and so can be
computed from the traces of $M$ and $M^2$. The rotation
vectors of  fixed points are shown in Figure~\ref{rotset} as circles and 
those of the period two points by x's. 

\subsection{context}
The study of the simplest version of this rotation vector on the circle
goes back to Poincar\'e. The general case on manifolds
for flows was initiated by Schwartzman and Fried
\cite{schwartzman, fried1}. The adaption of these methods to iterated
maps requires just minor modification. There is a large
literature on certain case, especially on the two-dimensional
annulus and torus. For computing the rotation set using a lifted
Markov partition with references see \cite{Bdamster}. The general
case of an observable on a subshift of finite type is in 
 Ziemian \cite{zie}. 

\section{Fixed and periodic points}\label{sect4}
To compute the asymptotic displacement in the general case of $A\not=I$
 we first work with fixed and
periodic points and in the next section deal with general orbits

\subsection{Displacements of fixed points and the Bowen-Franks group}
Assume that $x\in\fix(f)$ and so for lifts  $\tx, \tf$  there
is an $n\in\Z^b$ with $\tf(\tx) = \deltan  \tx$. Using \eqref{deck},
$\Delta(\tx, \tf) = \vn$. Thus for fixed points we can obtain
the translation without using the map $\iota$. From \eqref{depends} to
get a translation independent of the choice of lift of $x$
 we define the \textit{displacement} of $x$ as
\begin{equation}
D(x, \tf) = \Delta(\tx, \tf) + (A-I)\Z^b \in \frac{\Z^b}{(A-I)\Z^b} := \BF_1(f)
\end{equation}
where $\tx$ is any lift of $x$. 
The \textit{fixed point displacement set} of $f$ is
$$
D_1(f) = \{D(x,\tf): x\in\fix(f)\}\subset \BF_1(f)
$$
where we suppress the choice of lift $\tf$.
The $\BF$ stands for Bowen-Franks
and the group $\BF_1$ occurs in many contexts,
 see section~\ref{BFsect} below.
 Using the Smith normal
form  $\BF_1(f)$ is a finite group
of order $\det(A-I)$ when one is not an eigenvalue of $A$.

\subsection{example 2}\label{eg2}
Consider the homomorphism $\psi_2:F_2\raw F_2$  induced by 
$a\mapsto a a a b, b \mapsto b b b a$. Let $\phi_2:X\raw X$ be 
its tight map and lift $\tphi_2$ as in Example~\ref{eg1}.
After computing the image of a fundamental domain
 the rest of the action of $\tphi_2$ can be determined by the formula 
$\tphi_2(\tx + \vn) = \tphi_2(\tx) + A_2 \vn$ where 
\begin{equation*}
A_2 = (\phi_2)_* = \begin{pmatrix} 3 & 1\\1 & 3\end{pmatrix}
\end{equation*}
is the action on $H_1(X, \Z) \equiv \Z^2$.

 Let
$\Gamma = (A_2-I)\Z^2$  and  
so $\BF_1(A_2) = \Z^2/\Gamma \cong \Z_3$ using the Smith
normal form. In the fundamental
domain  there are lifts of $5$ distinct fixed points,
$(0, 0), (1/3, 0), (2/3, 0), (0, 1/3)$, and $(0,2/3)$. Computing, 
$\phi(1/3,0) = (1/3, 0) + (1,0)$ and so $D(1/3,0) = (1,0) + \Gamma$.
Similar computations yield $D(2/3,0) = (2,0) + \Gamma, 
D(0,1/3) = (0,1) + \Gamma$, and $D(0,2/3) = (0,2) + \Gamma$. Now
note that $
(1,0) - (0,1) = (A_2-I) (1,-1)^T \in (A_2-I)\Z^2$ and so 
$D(1,0) = D(0,1)$. Similarly, $D(2,0) = D(0,2)$. Also,
$(1,0) - (2,0) \not\in (A_2-I)\Z^2$ and so $D(1,0) \not= D(2,0)$
and similarly $D(0,1)\not = D(0,2)$. Thus among
the $5$ fixed points of $\phi$ there are three distinct displacements
and so $D_1(\phi) = \BF_1(A_2)$.

\subsection{Comparing fixed points}
Let us explore the implications of two fixed points
having the same displacement. So assume $\tf(\tx) = \delta_{\vm}\tx$
and $\tf(\tx') = \delta_{\vm'} \tx'$ with $\vm' - \vm = (A-I)\vn$ for some
$\vn\in\Z^b$. In this case, 
\begin{equation}
\tf(\delta_{\vn}\tx) = \delta_{A\vn}\tf(\tx) = \delta_{A\vn} \delta_{\vm}\tx 
= 
\delta_{\vm + (A-I)\vn}\delta_{\vn}\tx  = 
\delta_{\vm'}(\delta_{\vn}\tx).
\end{equation}
Thus $\delta_{\vn}\tx$ and $\tx'$ are both fixed by the 
lift $\delta_{-\vm'}\tf$.
We have proved the nontrivial half of the fact that 
$D(x) = D(x')$ if and only if there are lifts $\tx, \tx'$ and $\tf$
with $\tf(\tx) = \tx$ and $\tf(\tx') = \tx'$.

Exploring further, consider an arc $\tgamma:[0,1]\raw \tMA$ with
$\tgamma(0) = \tx$ and $\tgamma(1) = \tx'$ with $\tx$ and $\tx'$ both
fixed by the same lift $\tf$. The arc $\tf\circ\tgamma$ also goes
from $\tx$ to $\tx'$ and so $\tgamma^{-1}\# (\tf\circ\tgamma)$ is a closed
loop in $\tMA$ and so pushing down to $M$, $\gamma^{-1}\# (f\circ\gamma)$
is homologous to zero in $M$ and so $\gamma$ is homologous to 
$f\circ\gamma$ rel endpoints. Following the steps in the other direction
yields the converse and so we summarize

\begin{lemma}\label{first3}
For $x, x'\in \Fix(f)$ the following are equivalent.
\begin{itemize}
\item $D(x, \tf) = D(x',\tf)$ for some lift $\tf$.
\item There are lifts $\tx, \tx'$ and $\tf$ with $\tf(\tx) = \tx$
and $\tf(\tx') = \tx'$.
\item There is an arc $\gamma$ from $x$ to $x'$ so that 
$\gamma$ is homologous to $f\circ \gamma$ rel endpoints.
\end{itemize}
\end{lemma}

\subsection{example 2, continued}\label{eg2a}  
Recall that $\Delta(1/3, 0) = (1,0)$ and $\Delta(0,1/3) = (0,1)$
while $(1,0)^T - (0,1)^T = (A_2-I) (1, -1)^T$. Thus using the calculation above,
$(0, 1/3) + (1, -1)$ and $(1/3,0)$ should be both fixed by the lift
$\tphi_2 - (1,0)^T$ which is  also easy to verify directly. 

\subsection{periodic points} 
The generalization to periodic points is straightforward. First note
that not every lift of $f^k$ is of the form $\tf^k$ for some
lift $\tf$ of $f$. So fix  a lift $\tfk$ of $f^k$ and assume
$x\in\Fix(f^k)$ and for a lift $\tx$, let the translation be 
\begin{equation}
\Delta_k(\tx, \tfk) = \titfk(\tx) - \tiota(\tx).
\end{equation}
To get a displacement independent of the choice of lift, let
\begin{equation}
D_k(x, \tfk) = 
\Delta_k(\tx, \tfk) + (A^k -I)\Z^b \in \frac{\Z^b}{(A^k-I)\Z^b}
:= \BF_k(f).
\end{equation}
Note that at this point we are just requiring that $x\in\Fix(f^k)$
not that $x$ has least period $k$.
If $x$ is a least period $k$ point the displacements of other points
on its orbit are simply related by $D_k(f^i(x), \tfk) = A^i D_k(x,\tfk)$.
\begin{stand}\label{stand} Henceforth we assume that the linear transformation
$A$ has no eigenvalues which are roots of unity. This implies that
$A^k-I$ is invertible and 
that $\BF_k(A) = \frac{\Z^b}{(A^k-I)\Z^b}$ is a finite group
for all $k$.
\end{stand}

\subsection{example 2, continued}\label{eg2b}
For $\phi_2$,  
in contrast to the fixed point case, for $k>1$
there are far more possible displacements in $\BF_k(A_2)$ than
there are fixed points. The eigenvalues of $A_2$ are $2$ an $4$
and so the order of $\BF_k(A_2)$
is $\det(A_2^k-I) = 8^k - (2^k + 4^k) + 1$. So, for example, when
$k=2$ there are $45$ different possible displacements in $\BF_2(A_2)$.

On the other hand, using each circle in $X$ as a Markov box,
the dynamics of $\phi_2$ are described by the edge shift of
the matrix $A_2$. When $k$ odd, the vertex is counted twice
and when $k$ is even, the vertex is counted four times. Thus
$\#\, \Fix(\phi^k) = \trace(A_2^k) - N_1(k) = 
2^k + 4^k - N_1(k)$ where $N_1(k) = 1$ when $k$ is odd and
$N_1(k) = 3$ when $k$ is even. Thus, for example, when $k=2$
there are $18$ fixed points in contrast to the $45$ available
displacements.

\subsection{context} In the language of Nielsen Theory
the displacement as defined here is
the Abelianization of the twisted conjugacy class of a periodic
point and the displacement class is the Abelian Nielsen class.
 The order of $\BF_k(f)$ gives a lower bound on the number
of Nielsen classes of order $k$.
 
\section{The direct limit group $\BF_\infty$}\label{BFinfinity}
\subsection{Defining the direct limit}
As noted in the introduction the goal of this paper  is to build invariants
of a dynamical system which capture information about 
asymptotic dynamics using first homology. The next step is to collect 
together the displacements of all the periodic
points. The crucial, but elementary, observation is that if $x\in\Fix(f^i)$
and $i\vert j$ ($i$ divides $j$), then $x\in\Fix(f^j)$. The algebraic
analog of this fact allows us to build $\BF_\infty$, the direct limit
of the groups $\BF_i$. Since initially the consideration are strictly
 algebraic we work with a fixed square integer matrix $A$ which satisfies
the standing hypothesis~\ref{stand}.

Let $(\N, \prec)$ be the directed set of positive integers
 ordered by divisibility,
so $i\prec j$ if and only if $i\vert j$. When $i\prec j$ define
\begin{equation}
C_{ij} = A^{(\frac{j}{i} -1)i} + A^{(\frac{j}{i}-2)i} + \dots A^i + I
\end{equation}
and so $A^j-I = C_{ij} (A^i-I)$. Let $\Gamma_i = (A^i-I)\Z^b$ and 
for $i\prec j$ define
$\Upsilon_{ij}: \BF_i(A)\raw \BF_j(A)$ via
\begin{equation}
\vn + \Gamma_i \mapsto C_{ij}\vn + \Gamma_j
\end{equation}
It is easy to check that  $\Upsilon_{ij}$ is a well-defined monomorphism and 
further, when $i\prec j \prec k$ we have $C_{jk} C_{ij} = C_{ik}$. Thus 
\begin{equation}
\Upsilon_{jk} \Upsilon_{ij} = \Upsilon_{ik}
\end{equation}
and so $(\{B_i\}, \{\Upsilon_{ij}\}, (\N, \prec))$ is a directed system
of Abelian groups and we define
\begin{equation}
\BF_\infty(A) = \varinjlim \BF_i.
\end{equation}
We treat the direct limit as 
\begin{equation}
\BF_\infty(A) = (\sqcup \BF_i)/\sim
\end{equation}
where $g_i \in \BF_i$ and $g_j\in \BF_j$ are equivalent,
$g_i\sim g_j$, if and only if when $i\prec k$ and $j\prec k$
then $\Upsilon_{ik}(g_i) = \Upsilon_{jk}(g_j)$ in $\BF_k$. We 
denote equivalence classes as $[\cdot]$. Recalling that
elements of the $\BF_i$ are cosets using the 
definitions of the $\Upsilon_{ij}$ the product in $\BF_\infty$ is 
\begin{equation}
[\vn + \Gamma_i] \cdot [\vn' + \Gamma_j] =
[(C_{ik}\vn + C_{jk}\vn') + \Gamma_k]
\end{equation}
when $i,j \prec k$. 

\begin{remark}\label{remdirect}  When $i\vert j$  one can directly
 compare the displacement
of $x\in\Fix(f^i)$ in $\BF_i$ with the displacement of the same point 
 as an element of $\Fix(f^j)$.
For $x\in\Fix(f^i)$ we have 
$\titf^i(\tx) = \titx + \Delta_i(\tx)$   and so 
$\titf^{2i}(\tx) = \titx+ \Delta_i(\tx) + A \Delta_i(\tx)$.
Continuing if $i\vert j$ we have $\titf^{j}(\tx) = \titx +
C_{ij}\Delta_i(\tx)$ which is to say that 
$\Delta_j(\tx) = C_{ij} \Delta_i(\tx)$ in agreement
with the definition of $\Upsilon_{ij}$.
\end{remark}

Returning to the dynamics   of $f:M\raw M$ 
we collect together
all the displacement data of all periodic points in 
\begin{equation}
D_\infty(f) = 
\bigcup_{n\in\N} \bigcup_{x\in\Fix(f^n)} [D_n(x)]\subset\BF_\infty(A).
\end{equation}

\subsection{embedding $\BF_\infty(A)$ in $\T^b$}\label{intotorus}
In this subsection we give a concrete realization of $\BF_\infty(A)$
as a subgroup of the torus group $\T^b = \R^b/\Z^b$ and in the
next subsection discuss it dynamical significance.

Fix a matrix $A$ satisfying the standing hypothesis.
For all $i$, define $\Psi_i:\BF_i(A)\raw \T^b$ via
\begin{equation}
(\vn + \Gamma_i) \mapsto (A^i-I)^{-1}\vn + \Z^b.
\end{equation}
This is easily seen to be well-defined. Now define
$\Psi:\BF_\infty(A)\raw \T^b$ via
\begin{equation}
[\vn + \Gamma_i] \mapsto \Psi_i(\vn + \Gamma_i)
\end{equation}
To check it is well-defined, assume that $[\vn + \Gamma_i] =
[\vn' + \Gamma_j]$ in $\BF_\infty(A)$. This means 
that there is a $k$ with $i\vert k$ and $j\vert k$ and
a $\vm\in\Z^b$ with
\begin{equation}\label{e1}
C_{ik}\vn - C_{jk}\vn' = (A^k-I)\vm.
\end{equation}
Now $C_{ik} = (A^i-I)^{-1} (A^k - I)$ and so substituting
into \eqref{e1} yields
\begin{equation}\label{again}
(A^i-I)^{-1}\vn - (A^j-I)^{-1}\vn' = \vm\in\Z^b
\end{equation}
and so $\Psi([\vn + \Gamma_i]) = \Psi([\vn' + \Gamma_j])$
in $\T^b$.

To compute the kernel of $\Psi$, say that $\Psi([\vn + \Gamma_i]) = 0$
which is to say that $(A^i - I)^{-1}\vn = \vm \in \Z^b$. But then
$\vn = (A^i-I)\vm$ so $\vn\in\Gamma_i$ and thus $[\vn + \Gamma_i] = 0$
in $\BF_\infty(A)$ and $\Psi$ is a monomorphism.

\begin{remark} The subgroup $\Psi(BF_\infty(A))\subset\T^b$
is isomorphic to the direct limit group
$\BF_\infty(A) = \directlim \BF_i(A)$ but
it is worth noting that the topology of $\Psi(\BF_\infty(A))$ induced
as a subset of $\T^b$ does not correspond to the direct limit topology on
$\BF_\infty(A)$ which is discrete. 
\end{remark}

\subsection{Toral dynamics and $\BF_\infty{A}$}\label{toral}
Fix an integer matrix $A$ that satisfies the standing hypothesis and let 
$\Phi_A$ be the induced toral endomorphism $\Phi_A:\T^b\raw\T^b$. The
universal Abelian cover of $\T^b$ is also its universal cover, $\R^b$.
The linear transformation $A:\R^b\raw \R^b$ is a lift of $\Phi_A$ and
it will be the preferred lift when computing displacements.

We reverse our point of view and now instead of finding the displacement
of a periodic point we try to find a period $k$ point of $\Phi_A$ that has 
a given
displacement $\vn\in\R^b$. Thus 
we must solve $A^k \tx = \tx + \vn$ and so $\tx = (A^k-I)^{-1}\vn$. This
always exists because $A$ has no eigenvalues that are roots of unity.
Thus for all $\vn$ there is a $\tx$ with $\Delta_k(\tx, A^k) = \vn$
and so $D_k(x, A^k) = \vn + (A^k-I)\Z^b\in \BF_k(A)$. Thus the collection
of displacements of $\phi_A$ is all of $\BF_\infty(A)$ or
$\Psi(D_\infty(A)) = \Psi(BF_\infty(A))$. In words, 
 all displacement classes of linear toral endomorphisms are filled.

Next consider $\Fix(\Phi^k_A)$ inside the toral group $\T^b = \R^b/\T^b$. Define
$h:\Z^b\raw \T^b$ via $h(\vn) = (A^k-I)^{-1}\vn + \Z^b$ and so
$h(\Z^b) = \Fix(\Phi_A^k)$. Now $\vn \in \ker(h)$ if and only if 
$\vn = (A^k-I)\vm$ for some $\vm\in\Z^b$. Thus $\ker(h) = (A^k-I)\Z^b$
and so as groups, $\Fix(\Phi_A^k)\cong \BF_k(A)$. This also shows that
$\Fix(\Phi_A^k)$ is the image of $\BF_k(A)$ under
the monomorphism $\Psi_k$ of the previous section.

Now we consider the relationship of $\BF_\infty(A)$ with all the fixed
points of $\Phi_A$, $\cup_{i\in\N}\Fix(\Phi_A^i)$. Let $F_i = \Fix(\Phi^i_A)$
and note that if $x\in F_i$ and $i\vert j$ then $x\in F_j$. So
when $i\vert j$ define $\phi_{ij}:F_i\raw F_j$ by $\phi_{ij}(x) = x$.
Using the directed set $(\N, \prec)$ under divisibility, we have that
$(\{F_i\}, \{\phi_{ij}\}, (\N, \prec))$ is a directed system and
let
\begin{equation}
F_\infty = \directlim F_i.
\end{equation}
Since $\phi_{ij}$ is the identity and all the $F_i\subset\T^b$, then
$F_\infty = \cup F_i$, the collection of periodic points of $\Phi_A$,
is a subgroup of $\T^b$. We next see that this subgroup
is isomorphic to $\BF_\infty(A)$. 

The equation~\eqref{again} and the lines above 
it imply that $\phi_{ij}\circ \Psi_i = 
\Psi_j \circ \Upsilon_{ij}$ and so the map $\Psi$ in fact sends equivalence
classes in $\BF_\infty(A)= \directlim \BF_i$ to equivalence classes
 in $F_\infty = \directlim F_i$ and thus the image of
$\BF_\infty(A)$ under the isomorphism $\Psi$ is $F_\infty$,
the collection of periodic points of $\Phi_A$.

\begin{theorem}
Given the linear transformation $A$ satisfying the standing hypothesis,
 the collection of fixed points
of the toral endomorphism $\Phi_A$ is a subgroup of the torus
group $\T^b$ and can naturally be identified with the 
group $\BF_\infty(A)$ of possible displacements of $\Phi_A$. 
The identification
is simply that of a fixed point of $\Phi_A^i$ with its displacement.
\end{theorem}

As a consequence, considering a general map $f:M\raw M$ with
$H_1(M;\Z) = \Z^b$ and $f_* = A$, the collection of its
displacements is a subset of $\BF_\infty(A)$ and this may
also be viewed as a subset of the periodic points of the
total endomorphism $\Phi_A$. So the process of computing displacements
can be viewed as comparing the periodic points
 of $f$ with those of a toral endomorphism. We explore this
comparison more closely in the next section.

\section{Global shadowing}\label{sectGS}

\subsection{The comparison map $\alpha$}
The map $\tiota:\tMA\raw\R^b$ allows us to compare 
$\tf$ orbits and $A$ orbits.
In the bases, this corresponds to comparing orbits of $f$ and $\Phi_A$.
Specifically, we compare  the lift $\tx$
of a periodic point of $f$ to the point $\ty$ to which
is corresponds under the identification of the last section.
In symbols, if $\tfk(\tx) = \deltan\tx$ we know there is a
 point $\ty$ with $A^k\ty = \ty + \vn$, or $\Delta_k(\tx, \tfk) = 
\Delta_k(\ty, A^k)$. Let $\talpha(\tx) = \ty$ and so 
$\talpha$ maps lifts of periodic points of $f$ to lifts
of periodic point of $\Phi_A$. Recall that $\Delta_k(\deltan\tx, \tfk)
= \Delta_k(\tx, \tfk) + (A^k-I)\vn$ and similarly 
$\Delta_k(\deltan\tx, A)
= \Delta_k(\tx, A^k) + (A^k-I)\vn$. This implies that $\talpha$
descends to a map $\alpha$ from the periodic points of $f$ in $M$
to the periodic points of $\Phi_A$ in $\T^b$. The map $\alpha$
sends a period $k$ point for $f$ to the unique period $k$
point $\alpha(x)\in\T^b$ with $D_k(x, f^k) = D_k(\alpha(x), \Phi_A^k)$.

\subsection{global shadowing of periodic points}
Let $(x,f)$ and $(y, \Phi_A)$ be period $k$ points.
Assume that $\ty = \talpha(\tx)$ and so for some $\vn, k$
$\tf^k(\tx) = \deltan\tx$ and $A^k \ty = \ty + \vn$.
Given $m\in\N$ write it as $m = ik + j$ with $0\leq j < k$. 
Repeatedly using~\ref{fundamental}
 yields 
\begin{equation*}
\tf^{ik}(\tx) = \delta_{\vn + A^k\vn + \dots + A^{(i-1)k}\vn}\, \tx  
\end{equation*}
and so 
\begin{equation*}
\tf^m(\tx) =  \delta_{A^j(\vn + A^k\vn + \dots + A^{(i-1)k}\vn})\tf^{j}(\tx).
\end{equation*}
Thus
\begin{equation*}
\titf^m(\tx) = \titf^{j}(\tx) + A^j(\vn + A^k\vn + \dots + A^{(i-1)k}\vn)
\end{equation*}
and a similar formula for $A^m\ty$, and so 
\begin{equation}
\|\titf^m(\tx) - A^m \ty\| = \|\titf^j(\tx) - A^j\ty\| := K_j 
\end{equation}
The collection $\{K_j\}$ is finite so we conclude there
is a $K$ so that $\|\titf^m(\tx) - A^m \ty\| < K$ for all $m\in\N$.

\begin{definition}
Given a lift $\tf:\tMA\raw \tMA$
two orbits in the covers $(\tx, \tf)$ and $(\ty,A)$ 
are said to shadow if there
is a constant $K$ with 
\begin{equation}
\td(\titf^n(\tx),A^n(\ty)) \leq K\ \text{for all}\ n\in\N
\end{equation}
where $\td(\tx,\ty) = \|\tx-\ty\|$ in $\R^b$.
Two orbits in the bases $(x,f)$ and $(y,\Phi_A)$ are said  to
globally shadow if there are lifts $\tx, \ty, \tf$ 
that shadow.
\end{definition}
There are a number of uses of shadowing and in a general context
this form of shadowing would be called Abelian global shadowing.
The argument above shows
\begin{lemma}\label{pergs}
 Assume $f:M\raw M$ has $f_* = A$ acting on $H_1 \equiv \Z^b$
and $x$ and $y$ are period $k$ points of $f$ and $\Phi_A$, 
respectively. If $D_k(x, f^k) = D_k(y, \Phi_A^k)$,   or equivalently
  $\alpha(x) = y$,
then $x$ and $y$ globally shadow.
\end{lemma}

\subsection{example 2 continued}\label{eg2c}
Recall $D(1/3,0) = (1,0)$ and so $\talpha(1/3,0)
= (A-I)^{-1}(1,0) = (2/3, -1/3)$. Thus $o((1/3,0), \tf))$
shadows that of $o((2/3, -1/3), A)$ in $\R^2$. Similarly,
since $D(0,1/3) = (0,1)$, $o((0,1/3), \tf)$ shadows $o((-1/3, 2/3),A)$
in $\R^2$. Projecting to $X\subset \T^2$, note that 
$(2/3, -1/3) \equiv (2/3, 2/3) \equiv (-1/3, 2/3) \mod \Z^2$ and so
indicating projections by primes, $((2/3, 2/3)', \Phi_A)$,
$((1/3, 0)', f)$, and $((0,1/3)', f)$ all globally shadow.
Similarly, $((1/3, 1/3)', \Phi_A)$,
$((2/3, 0)', f)$, and $((0,2/3)', f)$ all globally shadow.

\section{$H_1$-expansion and general orbits}\label{sect7}
\subsection{$H_1$-expansion}
Since all periodic orbits of $f$ globally shadow a periodic orbit of
$\Phi_A$, the natural next step is to investigate general
orbits.  To make progress we need more control over
$A$.
\begin{definition}\label{Hone}
A linear transformation $A:\R^b\raw \R^b$ is called expanding if all 
its eigenvalues are outside the unit circle.
A map $f:M\raw M$ is $H_1$-expanding if $A=f_*$ is.
For expanding  $A$,
let $\lambda = \lambda(A)$ be such that $\rho(A^{-1}) < 1/\lambda < 1$
and  $\|\cdot \|$ be an adapted norm with $\|A^{-1}\|<1/\lambda$
 (see Proposition 1.2.2 in \cite{katokhasse}). Thus
$\|A^{-1}\vv\| < (1/\lambda)\|\vv\|$ and $\|A\vv\| > \lambda\|\vv\|$
for all $\vv\in\R^b$. Accompanying $A$ will always be the metric 
on $\R^b$ coming from the adapted norm, $\td(x,y) = \|x-y\|$.
\end{definition}

The crucial observation in shadowing is that
in the presence of expansion there is a uniform bound
on the distance between orbits that shadow. If
 orbits get too far away from each other, the expansion takes
over and their mutual distance goes to infinity. From
another angle, \eqref{desired} says that  on scales larger than 
$c(f) := \|\sigma\|$
 the map $f$ looks like $A$ and so is ``expanding on large scales''.
\begin{lemma}\label{gs1} Assume $f$ is $H_1$-expanding. The orbits
$(x, f)$ and $(y, \Phi_A)$ globally shadow if and only if they
globally shadow with constant $c(f)/(\lambda(A)-1) := \delta(f)$.
\end{lemma}

\begin{proof} Assume $(\tx, \tf)$ and $(\ty, A)$ globally shadow.
Since $\titf(\tx) = A \titx + \sigma(\tx)$
letting $\delta_k = \|\titf^k(\tx) - A^k\ty\|$,
\begin{equation}
\delta_{k+1} = \|A(\titf^k(\tx) - A^k\ty) + \sigma(\tf^k(\tx))\|
\geq \lambda \delta_k - c(f).
\end{equation}
 Consider the dynamical system
$z\mapsto \lambda z - c$. It has a fixed point at $\delta(f)$
and any point greater than this fixed point iterates to infinity.
Thus if some $\delta_k > \delta(f)$, then $\delta_k\raw\infty$
contrary to assumption. The converse is trivial.
\end{proof}

\subsection{global shadowing}
With the expansion hypothesis the phenomenon of global
shadowing holds not just for periodic orbits, but for
general orbits as well.
The proof is a simple
variant  of the standard hyperbolic shadowing lemma.
\begin{lemma}\label{bound}
Let $f:M\raw M$ be $H_1$-expanding with $\lambda$ and  $c(f)$
as defined above.  For $\tx_0\in\tMA$ there is a unique $\ty\in\R^b$
such that 
\begin{equation}\label{best}
\|\titf^i(\tx_0)- A^i(\ty)\|\leq \delta(f)
\end{equation}
for all $i\in\N$. Thus $(x,f)$ globally shadows $(y,\phi_A)$
with constant $\delta(f)$.
\end{lemma}
\begin{proof}
For each $i$, let $\tx_i = \titf^i(\tx_0)$ and so using \eqref{desired},
$\tx_i$ is a $c(f)$-pseudo-orbit for $A$. When 
 $B_i = B_\delta(x_i)$ for $z\in B_i$,
\begin{equation}
\td(z, Ax_{i-1}) \leq \td(z, x_i) + \td(x_i, A(x_{i-1})) 
= \td(z, \tx_i) + \td(\titf(\tx_{i-1}), A\tx_{i-1}) 
< \delta(f) + c(f).
\end{equation}
Thus $\td(A^{-1}(z), x_{i-1}) < \frac{\delta(f) + c(f)}{\lambda}
= \delta(f)$ and so $A^{-1}(\overline{B_i})\subset B_{i-1}$.

Therefore, $\{A^{-i}(B_i)\}$ is a nested family of topological
balls in $B_0$ of decreasing diameter so there is a unique
point $\ty\in\cap A^{-i}(\overline{B_i})$.
 Thus $A^{i}(\ty) \in B_i$
for all $i$ and so $o(\ty,A)$ shadows $o(\tx, \tf)$ with
constant $\delta$.

For the uniqueness, assume that both $\ty$ and $\ty'$ satisfy
\eqref{best}. Then $\|A^i(\ty) - A^i(\ty')\| < 2 \delta$
which implies $\ty=\ty'$ since $A$ is expanding.
\end{proof}

\subsection{the semiconjugacy}
\begin{definition} Define $\talpha:\tMA\raw\R^b$ by letting
$\talpha(\tx)$ be the unique point $\ty$ with 
$o(\tx,\tf)$ shadowing $o(\ty, A)$.
\end{definition}

\begin{theorem}\label{talpha}
When $f$ is $H_1$-expanding, the map  $\talpha:\tMA\raw\R^b$ 
is continuous, equivariant ($\talpha(\deltan\tx) = \talpha(\tx) + \vn$),
and satisfies $\talpha\circ\tf = A \circ \talpha$. It thus descends
to a continuous semiconjugacy $\alpha:M\raw \T^b$ of $f$ acting on $M$ to
$\Phi_A$ acting on the compact, invariant set $\alpha(M)\subset \T^b$.
Further, $\alpha$ is homotopic to $\iota$ and thus if $M = \T^b$,
$\alpha$ is onto.
\end{theorem}
\begin{proof}
We first prove continuity. If $\tx_n\raw\tx_0$ we must 
show that $\talpha(\tx_n)\raw \talpha(\tx_0)$.

Let $\ty_n = \talpha(\tx_n)$. By Lemma~\ref{bound},
$\td(\tiota(\tx_n), \ty_n) \leq
\delta(f)$ for all $n$. Thus for sufficiently large $n$, 
$\td(\ty_n, \tiota(\tx_0)) < 2\delta$. Thus passing to a subsequence
if necessary, there exists a $\ty_0$ with $\ty_n\raw\ty_0$.
For a fixed $k$
\begin{equation}
\td(\titf^k(\tx_0), A^k\ty_0) \leq
\td(\titf^k(\tx_0), \titf^k(\tx_n)) + \td(\titf^k(\tx_n), A^k\ty_n)
+ \td(A^k\ty_n, A^k\ty_0).
\end{equation}
Letting $n\raw\infty$ yields $\td(\titf^k(\tx_0), A^k\ty_0)\leq \delta(f)$.
This holds for all $k$, and so $\ty_0 = \talpha(\tx_0)$ as 
required.

For the equivariance, note that \eqref{fundamental} works for both $\titf$
and $A$ and the metric is equivariant  and so 
$\td(\titf^k(\deltan\tx), A^k(\talpha(\tx)+ \vn) = 
\td(\titf^k(\tx), A^k\talpha(\tx)) \leq \delta(f)$ for all $k$, and so
$\alpha(\deltan\tx) = \talpha(\tx) + \vn$ for all $\vn\in\Z^b$.

As for the commutativity, note that $\td(\titf^k(\tx), A^k\talpha(\tx)) \leq
\delta(f)$ for all $k$ says that 
\begin{equation}
\td(\titf^{k-1}(\tf(\tx)), A^{k-1}(A\talpha(\tx)) 
\leq \delta(f)
\end{equation}
 for all $k$, and thus $\talpha(\tf(\tx)) = A\talpha(\tx)$
for all $\tx$. The rest of the second sentence follows immediately.

For the last sentence, let $\gamma(\tx) = \talpha(\tx) - \titx$.  Then 
inside $\R^b$, for $t\in [0,1]$, $\talpha_t(\tx) = \tiota(\tx)
 + (1-t) \gamma(\tx)$
is a homotopy from $\talpha$ to $\tiota$ that descends to
a homotopy from $\alpha$ to $\iota$ in $\T^b$. Since $\iota$ is
a homeomorphism when $M=\T^b$ and $\T^b$
 is a closed manifold, 
$\alpha:\T^b\raw\T^b$ is onto. 
\end{proof}

\begin{remark}
The semiconjugacy $\alpha$ is usually not onto, except in the special
case noted in the theorem. It is sometimes injective, but this can be 
quite difficult to prove. In addition,
$\alpha$ is often of low regularity, for example H\"older
but nowhere locally  of bounded variation. See 
Sections~\ref{reg} and \ref{hirsch} below.
\end{remark}

\subsection{context}
The notion of global shadowing is due to Katok. It is central
to the proof of Handel's Theorem~\ref{minthm} below. The arguments
and results of this section are Abelian versions of Handel's
from \cite{handel}. J. Kwapisz pointed
out to the author that  
the standard pseudo-orbit shadowing argument
  applied to  $\Phi_A$ lifted to $\R^b$ would yield
a semiconjugacy equivalent to that of
Franks of the next section (Theorem~\ref{FS1}). The shadowing
argument used here is a simplification of Robinson's \cite{clark}.

\section{asymptotic displacement I}
Global shadowing provides on method of computing asymptotic
displacements of general orbits. In this section we work more directly
by taking a normalized limit of the motion after $n$ iterates.

  Iterating \eqref{desired} yields
\begin{equation}\label{itdesired}
\titf^n(\tx) = A^n\titx + A^{n-1}\sigma(x) +
 \dots + A \sigma(f^{n-2}(x)) +
\sigma(f^{n-1}(x)).
\end{equation}
This indicates that when $f$ is $H_1$-expanding 
most orbits run off to infinity in $\tMA$ at 
the rate of $\|A^n\|$. Normalizing by that  rate yields
\begin{equation}
A^{-n}\titf^n(\tx) = \titx + A^{-1}\sigma(x) + A^{-2}\sigma(f(x))
+ \dots + A^{-n}(f^{n-1}(x)) := \titx + \sigma_k(x).
\end{equation}
Since $A$ is expanding, using the adapted metric from Definition~\ref{Hone},
$\| A^{-k}\sigma(f^{k-1}(x))\| < c(f)/\lambda^k$,
and so by the Weierstrass M-test, the limit as $n\raw\infty$ converges
to a continuous function
\begin{equation}\label{series}
\tbeta(\tx) := \lim_{n\raw\infty} A^{-n}\titf^n(\tx) = \titx + 
\sum_{n=1}^\infty A^{-n}\sigma(f^{n-1}(x)) := \titx + \sigma_\infty(x)
\end{equation}
where both $\|\sigma_k(x)\|, \|\sigma_\infty(x)\| \leq \delta(f)$.
So in a precise sense, $\tbeta(\tx)$ is measuring a normalized
asymptotic displacement of $\tx$ under $\tbeta$. We summarize
the properties of $\beta$:
\begin{theorem}[Franks, Shub, Band]\label{FS1}
When $f$ is $H_1$-expanding the map  $\tbeta:\tMA\raw\R^b$ 
is continuous and if $f$ is $L$-Lipschitz, $\tbeta$ is $\nu$-H\"older
for any $0 < \nu < \min(\log(\lambda(A))/\log(L), 1)$. Further, 
$\tbeta$ is equivariant, $\tbeta(\tx + \vn) = \tbeta(\tx) + \vn$,
and satisfies $\tbeta\circ\tf = A \circ \tf$. It thus descends
to a continuous semiconjugacy $\beta:M\raw \T^b$ of $f$ acting on $M$ to
$\Phi_A$ acting on the compact, invariant set $\beta(M)\subset \T^b$.
Further, $o(\tx, \tf)$ globally shadows $o(\ty, A)$ if and
only if $\ty = \tbeta(\tx)$ and thus $\tbeta = \talpha$ from
Theorem~\ref{talpha}.
\end{theorem}
\begin{proof}
As already remarked, continuity of $\beta$ follows from the Weierstrass M-test.
To prove the H\"older assertion, first note that since $f$ is Lipschitz,
then by \eqref{desired}, so is $\sigma$. Since bounded Lipschitz
functions are $\nu$-H\"older for any $0<\nu<1$, for some constant $L'$, 
$\|\sigma(x) - \sigma(y)\| \leq L' \|x-y|^\nu$.  We therefore have
$$\|A^{-n}\sigma(f^{n-1}(x)) - A^{-n}\sigma(f^{n-1}(y))\|
\leq \frac{\|\sigma(f^{n-1}(x)) - \sigma(f^{n-1}(y))\|}{\lambda^n}
\leq \frac{L'(L^{n-1}\|x-y|)^\nu}{\lambda^n}$$
Thus
$$\|\sigma_\infty(x)-\sigma_\infty(y)\| <
\big(\frac{L'}{\lambda}\sum_{n=1}^\infty (\frac{L^\nu}{\lambda})^{n-1}\big)
\|x-y\|^\nu,
$$ where the sum converges since by the choice of $\nu$, $L^\nu/\lambda < 1$. 
So $\sigma_\infty$ and thus $\tbeta$ are $\nu$-H\"older.

Since $\sigma(\deltan\tx ) = \sigma(\tx)$,  equivariance follows. 
For the commutativity,
\begin{equation}
\tbeta\circ\titf = \lim_{n\raw\infty} A^{-n}\titf^{n+1}
=  \lim_{n\raw\infty} A (A^{-(n+1)}\titf^{n+1}) = A\tbeta
\end{equation}

To prove the shadowing assertion, recall that 
$\tbeta(\tx) = \titx + \sigma_\infty(x)$ where $\|\sigma_\infty\|
\leq \delta(f)$. Thus $\td(\tbeta(\tx), \titx)\leq
\delta(f)$ and so 
\begin{equation}
\delta(f) \geq 
\td(\tbeta\circ\titf^n(\tx), \titf^n(\tx))
= \td(A^n \tbeta(\tx), \titf^n(\tx))
\end{equation}
which is to say that $o(\tbeta(\tx), A)$  shadows 
$o(\tx, \tf)$. For the converse
assume that $o(\ty, A)$  shadows 
$o(\tx, \tf)$. We have just shown that 
$o(\tbeta(\tx), A)$  shadows $o(\tx, \tf)$.
Thus $\|A^i(\ty)-A^i(\tbeta(\tx))\| < 2 \delta$ for 
all $i$ and so  $\ty = \tbeta(\tx)$ since 
$A$ is expanding.
\end{proof}

\subsection{example 2, continued}\label{eg2d}
\begin{figure}[htbp]
\begin{center}
\includegraphics[width=0.8\textwidth]{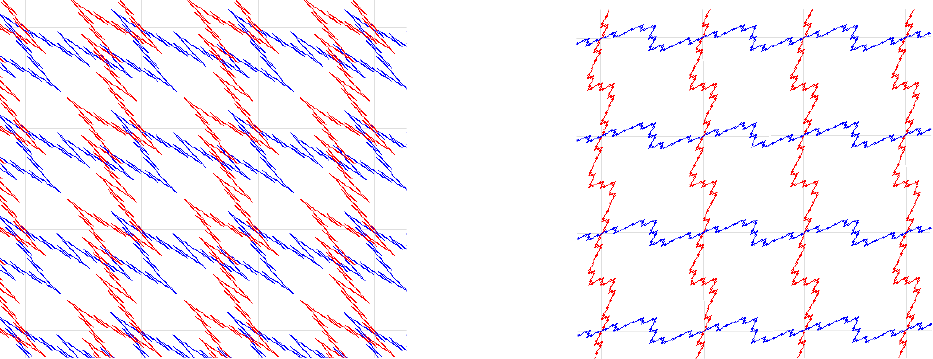}
\end{center}
\caption{Left: Approximation to $\tbeta(\tX)$ for the map $\phi_2$.
 Right: Approximation to $\tbeta(\tX)$ for the map $\phi_3$ of 
Example~\ref{Aexpand} below.
  }
\label{two_images}
\end{figure}

 Figure~\ref{two_images} (left)
 shows the image of $A_2^{-4} \tphi_2^4$ which is an
approximation to $\beta(\tX) = \alpha(\tX)$. See the Section~\ref{appendix}
for further details of the numerical method and  the
structure of the image.

\subsection{context} Franks proves Theorem~\ref{FS1} using a 
fixed point argument \cite{franks}. The summation
formula~\ref{series} is basically the $k^{th}$ iterate of his
operator acting on $f$. Shub uses a series in \cite{shub} and
shows the sum is H\"older when $f$ is Lipschitz. The equivalence
of Franks map to shadowing in the universal Abelian cover
is due to Band \cite{band}. The case of nontrivial Jordan
block is covered in \cite{bdeigen}. See Manning \cite{manning}
for some implications and amplifications of Franks Theorem.

\section{Asymptotic displacements II}\label{ADII}
We develop yet another natural way to compute asymptotic displacements.
In this method we take the limit of the approximate displacements for
finite orbit segments. It turns out that this is the same as
assigning a finite orbit segment a periodic point under 
$\Phi_A$ in $\T^b$ and then taking limits of these
 periodic points. Perhaps not surprisingly,
 the result agrees with the function $\talpha$
and $\tbeta$ of previous sections.

To approximate the orbit segments we partition $\tMA$ into
fundamental domains and use them for a coarse form of symbolic
dynamics. So let $\tM_0$ be a connected fundamental domain of $\tMA$ under the 
$\Z^b$ action by deck transformations and $\deltan\tM_{\vn} = \tM_0$
for $\vn\in\Z^b$ and so
$$\tMA = \sqcup \tM_{\vn}$$
as a disjoint union. In analogy with the situation on the
real line, for $\tx\in\tMA$, let $\floor{\tx} = \vn$ when
$\tx \in \tM_{\vn}$ and let $\{\tx\} = \titx - \floor{\tx} \in\R^b$.
 Thus for all $\tx$, $\titx = 
\floor{\tx} + \{\tx\}$ with $\|\{\tx\}\| \leq\diam{\iota M_0}$.
For each $\tx$ the approximate displacement after $k$ steps
is 
\begin{equation*}
d_k(\tx) = \floor{\tf^k(\tx)} - \floor{\tx}\in\Z^b
\end{equation*}
Thus $\titf^k(\tx) \approx \titx + d_k(\tx)$ with an error
uniformly bounded by twice the diameter of $\iota\tM_0$.
We imitate section~\ref{intotorus} and normalize by $(A^k-I)$ yielding
\begin{lemma}
When the function $f$ is $H_1$-expanding,
 $$\lim_{k\raw\infty} 
 (A^k-I)^{-1} d_k(\tx) = \tbeta(\tx).$$
\end{lemma}
\begin{proof}
For all $k$,
$$d_k(\tx) = \titf^k(\tx)-\{\tf^k(\tx)\} - (\titx - \{\tx\})
= (A^k-I)\titx + A^k\sigma_k(x) + \omega_k(\tx)$$
where $\omega_k(\tx) = -\{\tf^k(\tx)\} + \{\tx\}$ is bounded
by $2\diam{\tiota(M_0)}$.
Thus
\begin{equation}\label{approx}
(A^k-I)^{-1}d_k(\tx) = \titx + (A^k-I)^{-1}A^k\sigma_k(x) 
+ (A^k-I)^{-1}\omega_k(\tx) \raw \titx + \sigma_\infty(x) = \tbeta(\tx)
\end{equation}
as $k\raw\infty$ since $(A^k-I)^{-1}A^k = (I-A^{-k})^{-1}\raw I$
since $A$ is expanding.
\end{proof}

When $\tx$ is the lift of a period $k$ point then $d_k(\tx) = \Delta_k(\tx)$
and $\ty = (A^k-I)^{-1}d_k(\tx)$ is  the lift of a period $k$ point
of $\Phi_A$ and $(\tx, \tf)$ and $(\ty, A)$ shadow. This
means that downstairs, $\alpha(x) = y$, so informally we are using
$(y,\Phi_A)$ to approximate or model $(x, f)$. Thus for a general point
$\tx$ after $k$ iterates we have that $(A^k-I)^{-1}d_k(\tx)$ is the lift of
a period point of $\Phi_A$ that approximates $(\tx,\tf)$ for
$k$ iterates. Then \eqref{approx} takes the limit of these approximating
period points to get the probably not periodic point of $\Phi_A$ which
best approximates it. The lemma says that this point, as expected,
is $\alpha(x) = \beta(x)$.

\begin{remark} The definition of $d_k$ does not use the map $\tiota$.
The lemma above shows that the construction of the maps $\alpha=\beta$ does
not depend on the choice of the one-forms in the definition of $\iota$.
\end{remark}

\subsection{context} The method of just keeping track 
of the fundamental domain of a lifted orbit is common
in the  theory of rotation sets. See, for example,
Section 11 of \cite{Bdamster}.

\section{homeomorphisms}
For the sake of expositional simplicity thus far we have focused on
 $H_1$-expansion, but   most of the results are 
true, and were first proved, for homeomorphisms with 
a homological hyperbolicity assumption.
\begin{definition}
If $h:M\raw M$ is a homeomorphism and $(h_*)_1 = A$ has
all it eigenvalues off the unit circle, then $h$ is 
called $H_1$-hyperbolic
\end{definition}
It is essential that $h$ be a homeomorphism. Many of the results aren't true 
for a non-injective map even if the hyperbolicity
hypothesis on $(h_*)_1$ hold. For example, if $A$ is hyperbolic
but not invertible over $\Z$.

We briefly go through the previous sections and note what
changes are necessary. 
The material in 
 Sections~\ref{sect4} and  \ref{BFinfinity} just require that
$A$ have no eigenvalues which are roots of unity, so they
go through unchanged under $H_1$-hyperbolicity. Since we are
now working with a homeomorphism the definition of global
shadowing in Section~\ref{sectGS} needs to be extended to include
all iterates $n\in\Z$. For the material in Sections~\ref{sect7}
through Section~\ref{ADII} the strategy is to split the system
into stable and unstable directions and use the expansion
of $f$ and $f_*$ in the stable directions and using the
inverses, the expansion in the stable directions. This forces
several small changes, for example, the shadowing constant in 
Lemma~\ref{gs1} becomes $2 \delta(f)$ since both forward
and backward time need to be controlled. In
Theorem~\ref{FS1} when $f$ is $H_1$-hyperbolic,
 $\Phi_A$ is a hyperbolic
total automorphism, also called a linear Anosov map.

\section{Mixed spectrum}

\subsection{Single eigenvalues}
Understanding the influence of individual eigenvalues of $f_*$ requires
maps $\tMA\raw\R$. As in the  map $\tiota$, this requires
cohomology not homology. Somewhat informally, a cohomology class represented
by $H_1(M;\R)\raw \R$ is induced by a map $M\raw S^1$ which
lifts to $\tMA\raw\R$ satisfying certain conditions.

So assume now that $A^T = f^*:H^1(M;\R)\raw H^1(M;\R)$
has a real eigenvalue $\mu\not=0$ with eigenvector $\vv$. Let
$L_{\mu}:\R^b\raw \R$ be
$L_{\mu}(\tx) = \langle \vv, \tx \rangle$. Thus $L_{\mu}(A\tx) = 
\mu L_{\mu}$.
Using formula~\ref{itdesired}  and letting $\hsigma_\mu = L_\mu\circ \tsigma$
and $\sigma_\mu$ its projection to the base we have
\begin{equation}\label{single}
L_{\mu} \titf^n{\tx} = \mu^n  L_\mu \titx + \mu^{n-1}\sigma_\mu(x) +
 \dots + \mu \sigma_\mu(\tf^{n-2}(x)) +
\sigma_\mu(\tf^{n-1}(x)).
\end{equation}

\subsection{Single eigenvalues with $|\mu|>1$}
\subsubsection{The semiconjugacy}
Assume now that the eigenvalue of $A^T = f^*$ is real
with $\mu > 1$. The case of complex $|\mu|>1$ is easily dealt with.
 Dividing \eqref{single} by $\mu^n$ we have 
\begin{equation}\label{weier2}
 \lim_{n\raw\infty} \frac{L_\mu (\titf^n(\tx))}{\mu^n}
= L_\mu(\titx) +
 \sum_{i=0}^\infty \frac{\sigma_\mu (f^i(x))}{\mu^{i+1}}.
\end{equation}
This converges uniformly by the Weierstrass M-test, and 
let $\tbeta_\mu:\tM\raw \R$ be the continuous
function defined by this sum.
It follows directly  that 
\begin{equation*}
 \begin{CD}
 \tMA @>{\tf}>>\tMA\\
 @V{\tbeta_\mu}VV    @VV{\tbeta_\mu}V\\
 \R @>{\times\mu}>>\R%\\
 \end{CD}
 \end{equation*}
So  we have obtained a semiconjugacy
between $\tf$ acting on $\tMA$ and multiplication by
$\mu$ on $\R$.

\begin{theorem}[Franks, Shub]\label{FS2}
Assume that $f:M\raw M$ is a continuous map of the smooth, 
connected manifold
$M$ 
and $\mu\in\R$ or $\C$ is a simple eigenvalue of 
$f^*:H^1(M;\Z) \raw H^1(M;\Z)$
with $|\mu|> 1$.
 For each lift
$\tf:\tMA\raw \tMA$ of $f$ to the universal Abelian cover  $\tMA$,
there exists a unique  continuous map $\tbeta_\mu :\tMA\raw\R$ or $\C$  
 with 
\begin{equation*}
\tbeta_\mu \circ\tf = \mu \, \tbeta_\mu. 
\end{equation*}
If $f$ is $L$-Lipschitz, $\tbeta_u$ is H\"older with any constant
$0 < \nu < \log(|\mu|)/\log(L)$.
If $f$ is a \homeo\ we may consider $f\I$, and so
obtain $\tbeta_\mu$ for $0< |\mu| < 1$.
\end{theorem}

\begin{remark}\label{Jordan}
 If $A^T$ has a nontrivial Jordan block with eigenvalue
$\mu$  and generalized eigenbasis $A\vv_1 = \mu \vv_1,
A\vv_i = \mu\vv_i + \vv_{i-1}$ for $i= 2, \dots, n$ then one
obtains a family $\tbeta_1\circ\tf = \mu\tbeta_1, 
\tbeta_i\circ\tf = \mu\tbeta_i + \tbeta_{i-1}$ for $i= 2, \dots, n$.
\end{remark}

\subsubsection{Induced structures}
The semiconjugacy corresponding to a single
eigenvalue $\mu$ measures the asymptotic motion of 
orbits in  the covering space $\tMA$ in the direction of
the eigenvector. It does not, in general,
descend to the base but it does, however, induce two
structures that do descend to the base.

First, consider the level sets $ \tX_r = \tbeta_\mu\I(r)$.  
Since, $\tbeta_\mu \circ\tf = \mu \, \tbeta_\mu$,
we have $\tf(\tX_r) = \tX_{\mu r}$ and so the decomposition
of $\tMA$ into level sets of $\tbeta$ is $\tf$-invariant.
This decomposition descends to a $f$-invariant decomposition
$\{ X_r\}$ of $M$. In general the decomposition elements can
be quite wild unless $f$ has more structure, for example in the
pseudoAnsosov maps considered below.

Second, given a path $\tgamma:[0,1]\raw \tMA$, let 
$\tF_\mu(\tgamma) = \tbeta_\mu(\tgamma(1)) - \tbeta_\mu(\tgamma(0))$.
The function $\tF_\mu$ descends to a function $F_\mu$ which
assigns numbers to paths in $M$ in a way that is
homological. It is  a \textit{path cocycle}  which is defined as
 a continuous map 
$F:C([0,1], M)\raw \R$ ($\C$)  that is
additive $F_\omega(\gamma_1 \# \gamma_2) = 
F(\gamma_1) +  F(\gamma_2)$ and 
homology invariant; $\gamma_1$ homologous to 
$\gamma_2$ rel endpoints implies that $F(\gamma_1) =
 F( \gamma_2)$.

An path cocycle is said to
 represent the cohomology class
$c$,  if $F([\Gamma]) = L_c([\Gamma])$ for 
every closed curve $\Gamma$ with   homology class 
$[\Gamma]\in H_1(M;\R)$ and $L_c:  H_1(M;\R)\raw \R$ is
the linear functional representing the cohomology class $c$.
A path cocycle $F$ is acted on by the function
$f$ as $f^*(F)(\gamma) = F(f\circ\gamma)$.  
The prototypical example of a path cocycle is  the
assignment $\gamma\mapsto  \int_\gamma \omega$ for a closed
one-form $\omega$. A path cocycle is \textit{eigen} for
$f$ if $f^*(F) = \mu F$.

From Theorem~\ref{FS2} we get:
\begin{corollary}\label{cocycle}
 Assume $f:M\raw M$ is a continuous map and $\mu$ is an
eigenvalue $|\mu| > 1$
with eigen-class $c_\mu$ of the action of $f$ on first cohomology. 
There exists an eigen-path cocycle $F_\mu$ with $f^*(F_\mu) = \mu F$ and
$F_\mu$ represents the class $c_\mu$. 
\end{corollary}

\begin{remark}\label{cocycle2}
The eigen-path cocycle of the Corollary is closely related to the
invariant decomposition. Specifically, 
given an path cocycle $F$, define a  relation
on  $M$  by $z_1 \sim z_2$ if there exists an path $\gamma$
connecting $z_1$ and $z_2$ with $F(\gamma) = 0$.
This yields an equivalence relation on $M$ whose
classes define a decomposition of $M$.
For the eigen-path cocycle $F_\mu$, this is the
same as  the $f$-invariant decomposition $\{X_r\}$ defined
above.
\end{remark}

In the smooth case, given a closed one-form $\omega$
with $f^*\omega = \mu\omega$, integration along paths
gives an eigen-path cycle and the kernels of the one-form are tangent
to the invariant decomposition. This is a rare occurrence. Even
smooth maps will rarely have an eigen-one form, but rather an
eigen-cohomology class and so $f^*\omega = \mu\omega + d\chi$
for some smooth function $\chi$. The collection of path cocycles
``completes'' the space of one-forms
into a space of objects with lesser regularity
 where the action of $f$ on the  completion always has
an actual eigen-object, not just an eigen-equivalence class.

\subsection{Single eigenvalues $\mu = 1$}
If $A^T$ has an eigenvalue of $\mu = 1$ then \eqref{single} becomes
\begin{equation}\label{single2}
L_{\mu} \titf^n{\tx} = L_\mu \titx + \sigma_\mu(x) +
 \dots + \sigma_\mu(\tf^{n-2})(x) +
\sigma_\mu(\tf^{n-1}(x)).
\end{equation}
Since $\sigma_\mu$ is uniformly bounded, $L_{\mu} \titf^n{\tx}$
grows at most linearly.
Letting $\Delta^\mu(\tx) = L_\mu(\tf(\titx)) - L_\mu(\titx)$ we have 
as in Section~\ref{translation} that $\Delta^\mu(\tx) = \tsigma_\mu(\tx)$ 
and that $\Delta^\mu$ is independent of the choice
of lift of $x$. Thus we treat $\Delta^\mu(x) = \sigma_\mu(x)$ 
as a function on $M$.
Again as in Section~\ref{translation} it is natural to take the Birkhoff average
and define the rotation vector in the eigen-direction as
\begin{equation}\label{birk2}
\rho_\mu(x) = \lim_{n\raw\infty}\frac{1}{n} \sum_{i=0}^{n-1} \Delta_\mu(f^i(x))
= \lim_{n\raw\infty}\frac{L_\mu \titf^n(\tx) - L_\mu \titx}{n}.
\end{equation}
This limit  exists at  almost every point of an invariant
measure by Birkhoff's Pointwise Ergodic Theorem and measures
the asymptotic linear speed in the eigen-direction. 

\subsection{Additional cases}
The simplest case of mixed spectrum is when all the eigenvalues
of $A^T$ are either one or else off the unit circle. Then one
can use the methods above to compute rotation vectors or
semiconjugacies in the appropriate directions.

In Remark~\ref{Jordan} we discussed the results when $|\mu|> 1$ has
a nontrivial Jordan block. When $\mu=1$ we get a
$\rho_\mu$ in the direction of $\vv_1$ with $A^T\vv_1 = \vv_1$.
For other vectors in the generalized eigenbasis, there doesn't
seem to be much structure to exploit using the methods here.

If $A^T$ has eigenvalues that are roots of unity for some 
$k$,  they become eigenvalues of one for $(f^k)^* = (A^T)^k$ and we can
compute rotation vectors in those directions for $f^k$.

Since $A^T$ is the matrix of the action of $f$ on $H^1$, its characteristic
polynomial is monic with integer coefficients. Thus a classic theorem
of Kronecker says that if all the eigenvalues are on the unit circle
they are roots of unity.
But it can happen that $A^T$ has eigenvalues on the unit circle
which are not roots of unity. Borrowing a remark from \cite{salem},
if the characteristic polynomial is Salem by definition it is irreducible and
has one eigenvalue inside, another outside and the rest on the
unit circle. Since the polynomial is irreducible the roots on
the circle can't be roots of unity. In the theory we develop here,
it is not clear how these eigenvalues on the circle can be utilized
(\textit{cf.} \cite{fried3}).

A final case where the theory here fails is for  eigenvalues
inside the unit circle when $f$ is not a homeomorphism. 

\subsection{examples}\label{mixed}
As a first example, say $f:\T^2\raw \T^2$ acts on $H_1$ with repeated
one as an eigenvalue and a nontrivial Jordan block (or shear) 
\begin{equation}\label{shear}
A =\begin{pmatrix} 1 & 1\\0 &1\end{pmatrix}.
\end{equation}
Then $A^T$ has eigenvector $\vv = [0 \, 1]^T$ and so we can
compute rotation vectors of lifts
in the vertical direction in $\tMA = \R^2$, but in the horizontal
direction in the absence of expansion or contraction there
is nothing to exploit.

As a second example,
let $D_n$ be the closed two-dimensional disk with $n$ interior open
disks removed. We assume that $f:D_n\raw D_n$ is a homeomorphism
that keeps the outer boundary fixed. If $f$ fixes all the
disks then $A=I$ and we can defined a $n$-dimensional rotation
vector. Each dimension measures the asymptotic rate of rotation
or linking about one of the removed disks. In general, $f$ will permute
the disks in which case $A^T$ has a single eigenvalue of one with
eigenvector $\vv = [1 1 \dots 1]^T$. Thus there is a single dimension
of rotation vector which measures the asymptotic linking
of an orbit around all the disks. It is usual in this situation 
to use $\vv$ to create a $\Z$-covering space to which all $f$ lift.
Specifically, letting $L_\vv:H_1(D_n;\Z)\raw\Z$ be given by
$L_\vv(\vw) = \langle \vv, \vw \rangle$ and let $p:\tD_n\raw D_n$ be the
corresponding $\Z$-cover (it has $p_*(H_1(\tD_n;\Z)) = \ker L_\vv)$.
The one-dimensional rotation vector then measures the speed of orbits
moving up and down in this cover.

\subsection{Context} Fathi develops path cocycles and connects them
to a proper cohomology theory \cite{fathi}. Shub does Franks'
semiconjugacy for individual eigen-values using
Alexander cocycles and uses them to define the
invariant decomposition \cite{shub}. The contents of this section
are mainly taken from \cite{bdeigen}.  A fair number of authors have
considered the one-dimensional rotation set in the shear isotopy class
of the first example \cite{doeff1, doeff2, tali, sal1, salv}. 
The second example  is the context
for the Burau representation in which one records the action
of the map on the $\Z$-cover using homology with
 $\Z[t,t\I]$ coefficients with $t$ the deck transformation
(see, eg.\ \cite{bb}).

\section{regularity}\label{reg}
The resemblance of \eqref{weier2} to the definition of a
Weierstrass nowhere differential function is a hint that the
semiconjugacies $\beta$ and $\beta_\mu$ will, in general, have low
regularity. This is also expected from ``Conservation of Difficulty":
their existence is quite general and follows from simple,
robust hypothesis on the action on $H_1$ and the proof is just
one or two lines, so there has to be a catch somewhere.

The general heuristic is that the semiconjugacies often go from
a system of higher entropy to lower entropy. In practise this
means from more stretching to less stretching. Thus the semiconjugacy
must fold. In the cases studied, the image dynamics is transitive
and so the folds must take place on a dense set so one expects
$\beta$ in  these cases to not be locally of bounded variation
or differentiable.

\subsection{the circle}
The theory on the circle is simplified by using the monotone-light
decomposition of the semiconjugacy. 
A theorem  by Eilenberg and by Whyburn says that if $f:X\raw X$
is a continuous function of a compact metric space then there is a space
$Y$ and continuous maps $m:X\raw Y$ and $\ell:Y\raw X$ with $m$ monotone
(point preimages are connected) and $m$ light (point preimages are
completely disconnected) and $f = \ell\circ m$. In the case $X=S^1$
connected components of
point preimages are intervals and so the intermediate space $Y$ is also
the circle.

 If $f:S^1\raw S^1$ has degree $d>1$ and so $f_* = d$, by Theorem~\ref{FS1}
 there is a semiconjugacy $\beta:S^1\raw S^1 $ between $f$ and
 the map $d:\theta
\mapsto d\theta$. Since we are on the circle the semiconjugacy descends 
to a light semiconjugacy on the intermediate space, also the circle.
 Thus we lose
very little information by assuming that the given degree $d$ circle
map is light semiconjugated to $d$. In this case one can prove:
\begin{theorem}
Assume $f$ is a continuous, degree-$d$ circle map with a \mc{light}
semiconjugacy $\beta$. The following are equivalent:

\begin{itemize}

\item[(a)] The map $\tf$ is \mc{not} injective,
 
\item[(b)] The map $\beta$ is \mc{not} injective,
 
\item[(c)] There exists a full measure, dense, $G_\delta$-set
$\Lambda\subset S^1$
so that $\theta\in \Lambda$ implies that $\beta\Inv(\theta)$ is 
\mycol{completely disconnected and uncountable and thus contains a Cantor set}.

\item[(d)] The topological entropy of $f$ satisfies \mc{$\htop(f) > \log(d)$},

\item[(e)] For \mycol{all nontrivial intervals} $J \subset S^1$,
the map $\beta\vert_J$ is not of bounded variation.
\end{itemize}
\end{theorem}

\subsection{pseudoAnosov maps} For this section
we restrict to compact surfaces. A \textit{measured foliation} $\cF$
on $M$ is a one-dimensional foliation with a finite number of well-behaved
singularities called prongs and a holonomy invariant measure
on transverse arcs. A homeomorphism $\phi:M\raw M$ is \pA\ if there
are a pair of transverse measured foliations $\cF^u, \cF^s$
 and a number $\lambda>1$, called the dilatation, so that 
$\phi(\cF^s) = \frac{1}{\lambda} \cF^s$ and 
$\phi(\cF^u) = {\lambda} \cF^u$. This means that arcs in the stable
foliation $\cF^s$ are contracted by $1/\lambda$ and arcs in the unstable
foliation $\cF^u$ are expanded by $\lambda$. 

Assume now that $\phi$ is a \pA\ and the action $A^T = \phi^*$ on
$H^1$ has an unstable eigenvalue $\mu>1$. Let $F_\mu$
be the eigen-arc-cocycle give by Corollary~\ref{cocycle}.
If $\gamma$ is an arc wholly contained
in a  leaf of the stable invariant foliation then on
hand, $\diam(\phi^n(\gamma))\raw 0$ as
$n\raw\infty$
 on the other, if $F_\mu(\gamma) \not= 0$,
$|F_\mu(f^n(\gamma))| = |\mu^n F_\mu(\gamma)| \raw
\infty$.
 Thus of necessity $F_\mu(\gamma) = 0$ and so using  Remark~\ref{cocycle2},
 each decomposition element $X_r$
is composed of the union of leaves of the stable
foliation. That remark  also implies that $F_\mu$ is
holonomy invariant: if $\gamma$ is slid  to $\gamma'$ keeping
the endpoints on the same leaf, $F_\mu(\gamma) = F_\mu(\gamma')$.

What kind of transverse structure does $F_\mu$ define?
As a representative of an arc transverse to the
stable foliation, parameterize a segment of unstable leaf by its
arc length  $s$ and let $H(s) = F_\mu([0, s])$ 
($H$ is like a cumulative distribution function).  Now
$\phi$ stretches arclength by $\lambda$  and
since $\phi_* F_\mu = \mu F_\mu$,  the image of $H$ gets 
stretched by $\mu$.
 Thus 
\begin{equation}\label{scale}
H(\lambda s) = \mu H(s),
\end{equation}
Since every leaf in dense, this local
scaling of $H$ has to happen at a dense set of points $s$.

If $|\mu| < \lambda$, the everywhere local scaling given by \eqref{scale}
can only happen if $H$ \mc{folds up everywhere}.
  With some work this implies that while $H$ is   
 $\nu = \frac{\log(|\mu|)}{\log(\lambda)}$-H\"older, it is not H\"older for
any larger $\nu$, is
nowhere differentiable,  nowhere locally of
bounded variation, and nowhere locally injective.
This in turn implies that $F_\mu$ cannot
be extended from lengths of arcs to a measure.

However, since $H$ is $\nu :=\frac{\log(|\mu|)}{\log(\lambda)}$-H\"older,
one can make sense of 
$$
\Phi(\sigma) = \int \sigma \;dH,
$$
when $\sigma$ is $\nu_1$-H\"older with $\nu_1 > 1-\nu$.
Thus $F_\mu$ can be used to define a \mc{holonomy
invariant H\"older distribution} on transverse arcs to
the stable foliation.
 When $\mu = \lambda$, the function  $H$ 
is monotone and $\Phi$ is then a linear functional
on continuous functions on transverse arcs, so it is
just the unique, holonomy invariant transverse measure
to the foliation.

\begin{theorem}
Assume $\phi$ is a  \pA\ \homeo\ on the surface $M$ with
$\phi^*$ the action on $H^1(M;\Z)$. 
\begin{itemize}
\item For each unstable eigenvalue
$\mu$ of $\phi^*$ there exists an \mc{eigen-path cocycle} $F_\mu$
which is \mc{holonomy invariant} to the stable
foliation $\cF^s$. 
\item $F_\mu$ can be used to   
 define a transverse, holonomy
invariant  H\"older distribution  on $\cF^s$.
\item  This distribution is a measure if and only if
 $\mu = \lambda$, the dilation, (which can only
happen when $\phi$ has oriented foliations).
 \end{itemize} 

Similar results hold for each
 stable eigenvalue of $\phi_*$  
and  the unstable foliation, $\cF^u$.
\end{theorem}

\subsection{example}\label{twin}
Let $M$ be the genus two closed surface
and $\psi:M\raw M$ is a \pA\ map with characteristic 
polynomial of $\psi_*$ acting on
$H_1(M;\Z)$ which splits over the integers 
into a pair of irreducible quadratic factors with roots
$0 < \lambda\Inv < \mu\Inv < 1 < \mu <\lambda$ (recall that
$\psi_*$ is symplectic).
The eigenvalues/vectors yield \mc{four semi-conjugacies}
$\tbeta_\lambda, \tbeta_{\lambda\Inv}, 
\tbeta_\mu,$ and $\tbeta_{\mu\Inv}$.
Fathi (1988) shows that the Franks semiconjugacy into
$\T^4$ splits and descends into paired maps
$\beta_1 := (\beta_\lambda, \beta_{\lambda\Inv})$ and
 $\beta_2 := (\beta_\mu, \beta_{\mu\Inv})$, each a semiconjugacy
onto a linear, two-dimensional toral automorphism.

The semiconjugacy $\beta_1$ is a branched cover by \cite{franksrykken} and
so is locally a diffeomorphism at all but finitely many points and 
point inverses are finite sets. On the other hand, using the results
above  the semiconjugacy
$\beta_2$ is H\"older exponent
$\nu = \log(\mu)/\log(\lambda)$, but no larger  $\nu$. It is
nowhere differentiable and nowhere locally injective or B.V.
Typical point inverses are Cantor sets.

\subsection{context}
The circle results are from \cite{doubling} and the rest of the section
from \cite{bdeigen}. For an introduction to \pA\ maps
see \cite{Thurston, FLP, primer}. Fathi give a lower bound
for the capacity of  image of
 the semiconjugacy
described in the next section \cite{fathi}.

\section{Hirsch's problem and Fathi's program}\label{hirsch}
In 1970 Hirsch initiated the investigation into what kind of compact sets
can be invariant under a hyperbolic linear toral automorphism,
i.e. a linear Anosov map. There is a good amount of literature around
the problem and a nice summary is contained in \cite{fathi}. In short, 
if the invariant set is not a subtori, it must be fractal in various
 precise senses, for example,
contain no rectifiable arcs \cite{Mane}. One of the question left
open is whether the invariant set can be a non-toral ``very crinkled''
 submanifold?

In \cite{fathi} Fathi gave a strategy for attacking this question.
In its simplest version one starts with a \pA\ map $\phi$ on
a closed surface $\Sigma$ of genus $g>1$ which is $H_1$-hyperbolic. By Franks
Theorem~\ref{FS1}, there is a map $\beta:\Sigma\raw \T^{2g}$ whose image
is invariant under the action of $\Phi_A$ where $A = \phi_*$ on first
on homology. The question then reduces to
the injectivity of $\beta$.

To put this question into the context of this paper we give a general
definition which will also be used later in the paper.
\begin{definition}\label{coverexpand}
Let $f:M\raw M$ and $\tM$ be a covering space to which $f$ lifts. Say 
that $f$ is $\tM$-expanding if there is an equivariant metric $\td$ 
on $\tM$ which induces the manifold topology and a constant $\lambda>1$
with 
$\td(\tf(\tx), \tf(\ty)) \geq \lambda \td(\tx, \ty)
$ for all $\tx, \ty\in\tM$. 
\end{definition}
Note that since the metric is equivariant, the definitions are
independent of the choice of lift of the function.
For example, if $A$ is an expanding matrix then the toral
endomorphism $\Phi_A:\T^b\raw\T^b$
will be expanding in both the universal cover and the universal
Abelian cover (these are the same).  PseudoAnosov homemorphisms 
are  always hyperbolic in the universal cover (see section~\ref{GP}
 below). The question of whether \pA\ homeomorphisms are ever hyperbolic
in the universal Abelian cover, i.e.\ are 
 $\tMA$-hyperbolic, turns out to be the main issue in Fathi's program.

In general, when $f$ is $H_1$-expanding the next proposition says
that the injectivity of the semiconjugacy
$\beta$ is equivalent to $f$ being expanding in 
the universal Abelian cover, i.e. being $\tMA$-expanding.
Equation~\ref{desired}
 says that $H_1$-expansion implies large scale or coarse
expansion in $\tMA$. The injectivity of $\beta$ requires more,
namely, metric expansion in $\tMA$ at all scales. 
\begin{proposition}
Assume $f$ is $H_1$-expanding. The following are equivalent
\begin{itemize}
\item[(a)] The semiconjugacy $\beta:M\raw\T^b$ is
injective.
\item[(b)] $f$ is $\tMA$-expanding
\item[(c)] If $\td(\tf^i(\tx), \tf^i(\tx')) < K $ for some
$K$ and all $i$, then $\tx  = \tx'$, or informally,
no two points of  $\tMA$ shadow under $\tf$.
\end{itemize}
An analogous result is true for $H_1$-hyperbolic homeomorphisms
(see Remark~\ref{tMuhyper}). 
\end{proposition}
\begin{proof}
After observing  that $\tbeta = \tiota + \sigma_\infty$ with
$\sigma_\infty$ bounded implies that $\td(\tf^i(\tx), \tf^i(\tx'))
\raw \infty$ if and only if $\|\titf^i(\tx), \titf^i(\tx')\|\raw\infty$, 
the equivalence of (a) and (c) follows from Theorem~\ref{FS1}.
 Condition (b) easily
implies condition (c). We know $\beta$ is continuous and so assuming
(a) it is a homeomorphism onto its image since $M$ is compact.
Thus $\tbeta$ is a homeomorphism and it is easy to
check that $\td(\tx, \ty) = \|\tbeta(\tx) - \tbeta(\ty)\|$ 
is the required metric with $\lambda$ as in Definition~\ref{Hone}.
\end{proof}

\subsection{example 2 continued and example 3}\label{Aexpand}
The map $\phi_2:X\raw X$ from Example~\ref{eg2}
  is $\tXu$-expanding as will be shown
in Example~\ref{uexpand}
 below. However, in Example~\ref{eg2c} it 
was shown that there are pairs of points
that globally shadow in $\tXA$ and so the semiconjgacy into
$\T^2$ is not injective, as is seen in Figure~\ref{two_images}(left).

As an example of a $\tXA$-expanding map on $X$,
let $\psi_3$ be generated by $a\mapsto a a a b a a a,
b\mapsto b b b a b b b$ and $\phi_3:X\raw X$ be the tight map on a wedge of
circles.  Using Lemma~\ref{gs1}, since $c(\phi_3) < 1$ and $A$
in this case has eigenvalues $5$ and $7$,
  pairs of 
points in $\tXA$ can only globally
shadow and thus have the same $\beta$ image only if they globally shadow
with constant $\leq 1/4$. Thus it suffices to check pairs of points
within $1/4$ of points on a chosen fundamental domain. This is easily
done showing that $\beta$ is injective. 
See Figure~\ref{two_images}(right). One
can also show directly by checking cases that $\phi_3$ is
$\tXA$-expanding.

\subsection{Context}
Fathi \cite{fathi} showed that $\beta$ is injective
 in neighborhoods on non-singular 
points of the invariant foliations. Barge and Kwapisz \cite{bargekwap}
 showed $\beta$ is injective
on an open, full measure set in $\Sigma$. There are two cases
in which $\beta$ has been shown to be not injective: when
the  \pA\ map has just one singularity (Band \cite{band})
and when the surface is genus two and the invariant foliations
are orientable  (Bouwman and  Kwapisz \cite{BK}).

\section{The Bowen-Franks group and the suspension flow}\label{BFsect}
We have thus far encountered the Bowen-Franks group 
$\BF_1(A) = \Z^b/(A-I)\Z^b$ as the collection of all
possible displacements of fixed points of an $f:M\raw M$
with $A = f_*:H_1(M)\raw H_1(M)$ as well as the group of fixed points of
the linear toral endomorphism $\Phi_A$. The maps $\alpha= \beta$
send a fixed point of $f$ with given displacement to the fixed point
of $\Phi_A$ with the same displacement.
There are at least three others situations in which $\BF_1(A)$ arises.
 
Recall that given $f:M\raw M$ its mapping torus is 
$T_f = (M\times [0,1])/\sim$ with $(x, 1) \sim (f(x),0)$.
Its first homology is 
$H_1(T_f;\Z) =  \BF_1(A) \oplus \Z$. The suspension flow $\varphi_t^f$
 on $T_f$ is the projection of of the vertical flow on $M\times [0,1]$
to $T_f$. Periodic points of $f$ correspond to closed orbits of
$\varphi_t^f$. The $\Z$ component of $H_1(T_f)$ keeps track of the
period. For a fixed point, the first component is the displacement
of the fixed point.
\begin{lemma}
Two fixed points $x,x'\in\fix(f)$ have $D(x,f) = D(x', f)$ in
$\BF_1(A)$ if
and only if the corresponding closed orbits of $\varphi^f_t$ are
 homologous in $T_f$.
\end{lemma}

Recalling from Section~\ref{toral} that for a toral endomorphism, 
all the displacement
classes are full. Thus for its suspension flow every 
homology class with second coordinate equal to one contains a closed orbit
that goes once around, i.e.
corresponds to a fixed point. 
When $f$ is  $H_1$- expanding or hyperbolic 
the semiconjuagcy $\beta:M\raw\T^b$ induces a map
 $\beta':T_f\raw T_{\Phi_A}$ with the semi-flow on $T_f$ pushed
to one on  $\beta'(T_f)$ making it an invariant set of $\varphi^{\Phi_A}_t$.
A closed orbit that goes once around under
$\varphi^f_t$ is thus pushed to one of 
 $\varphi^{\Phi_A}_t$ and thus shares its homology
class. Therefore the semiconjugacies can be viewed as picking
out a subset of the homology classes in $T_{\Phi_A}$.

Similar considerations of course hold for points in $\fix(f^k)$ using
the mapping torus $T_{f^k}$ of $f^k$. A natural question is how the
higher period orbits lie in $T_f$, the mapping torus of $f$ itself.
Are they distinguished by homology? In general, no, since 
$H_1(M_f)$ can be quite ``small'. As a simple example,
consider the doubling map $d:S^1\raw S^1$ via $d(\theta) = 2\theta$.
In this case $H_1(M_d) = \Z$ and thus homology only sees the period.
But the number of periodic points of period $k$ grows exponentially
and they are all in different  displacement classes since every lift
of $d^k$ has exactly one fixed point. 
The situation for periodic orbits and homotopy in $T_f$ is quite
different; $\pi_1(T_f)$ does distinguish different Nielsen classes,
see Section~\ref{Nsusp}.

The group $\BF_1(A)$ also arises as the deck group of  the largest
cover to which $f$ lifts and commutes with all deck transformations
\cite{fried2}. Thus in this cover, every lift of $x$ is fixed
by the same deck transformation, namely, $D(x,f)\in\BF_1(A)$.
Finally, as it originally appeared,
 $\BF_1(A)$ is an invariant for flow equivalence
of the subshift of finite type with transition matrix $A$.

\subsection{example}\label{BFeg}
For a fixed $A$ the collection of groups $\BF_k(A)$ are 
somewhat mysterious (to the
author at least) and it is not clear how the algebra connects
to the dynamics in a meaningful way. For the matrix $A_2$ from
Example~\ref{eg2} the family is very regular,
$$\BF_k(A_2) \cong \Z_{2^{k}-1} \oplus \Z_{ 4^k -1}.$$
But other examples are more complicated, for example, for
$$
A_4 = \begin{pmatrix}2 & 1 \\ 1 & 1\end{pmatrix}
$$
$\BF_2(A_4) = \Z_5, \BF_3(A_4) = \Z_4 \oplus \Z_4, 
\BF_4(A_4) = \Z_3 \oplus \Z_{15}, \BF_5(A_4) = \Z_{11} \oplus \Z_{11},
\dots$.

\subsection{context}
The papers of Bowen and Franks \cite{BF1,BF2} were mainly
about $\BF_1(A)$'s role as a flow equivalence invariant when
$A$ is the transition matrix for a subshift of finite type. However,
it is clear from the papers that they also knew the connection
of their group to the homology of a mapping torus and to the
periodic points of the toral endomorphism $\Phi_A$.

\section{Homotopy and Persistence (Nielsen Theory)}
Thus far we have used first homology to measure the progress of
orbits and have therefore been focusing on displacements in the
universal Abelian cover $\tMa$. In this section we measure
progress of orbits  using the fundamental group and therefore
work with displacements in the universal cover $\tMu$.   These
considerations belong to Nielsen Theory. Most expositions follow
the history and develop Nielsen Theory first and the contents
of the paper thus far are considered later as the Abelianization.

 For simplicity of
exposition, assume $f$ has a fixed point $x_0$ which we use as the
basepoint for the fundamental group. We write
$\pi = \pi_1(M;x_0)$ and let $f_\#:\pi\raw\pi$ be the induced map.
Fix an  identification of $\pi$ with the deck group
of $\tMu$. If $\tf:\tMu\raw\tMu$ is
a lift of a continuous $f:M\raw M$ then 
\begin{equation}\label{withdeck}
\tf(\delta \tx) = f_\#(\delta) \tf(\tx)
\end{equation}
 for all $\delta\in\pi_1$.

\subsection{fixed points}
For $x\in\Fix(f)$ with a lift $\tx$, and a fixed lift $\tf$, if
$\tf(\tx) = \delta \tx$ then the $\tMu$-translation is
$\Delta_u(\tx, \tf) = \delta$. Note that if we choose a different
lift of $\tx$, namely $\gamma\tx$ then 
$$
\tf(\gamma\tx) = f_\#(\gamma) \tf(\tx) = f_\#(\gamma) \delta\tx
= f_\#(\gamma) \delta \gamma\I (\gamma \tx)$$
which is to say that $\Delta_u(\gamma \tx) = f_\#(\gamma)
 \Delta_u(\tx) \gamma\I$.
Thus to get a well defined displacement we define for $x\in M$
$$
D_u(x, \tf) = \{f_\#(\gamma) \Delta_u(\tx) \gamma\I : \gamma\in\pi_1\ 
\text{and}\ \tx \ \text{is any lift of}\ x\}.
$$
In Nielsen Theory  the $\tMu$-displacement is called
the \textit{twisted conjugacy} or  \textit{Reidemeister coordinate}.
 
As with the $\tMa$-displacement in Lemma~\ref{first3},
there are three equivalent ways of looking at $\tMu$-displacement.
\begin{lemma} Assume $x,y\in\Fix(f)$
 The following are equivalent
\begin{itemize}
\item $D_u(x, \tf) = D_u(y, \tf)$ for some lift $\tf$.
\item There exists an arc $\gamma$ connecting $x$ to $y$ with
$f\circ\gamma $ homotopic to $\gamma$ rel endpoints.
\item There exist lifts $\tx, \ty, \tf$ to $\tMu$  with $\tf(\tx) = \tx$
and $\tf(\ty) = \ty$.
\end{itemize}
\end{lemma}

The collection of fixed points of $f$ that have the same $\tMu$-displacement
is called the $\tMu$\textit{-displacement class} or 
the \textit{Nielsen class} of $x$. We adapt the latter classical terminology
here. The Nielsen class of $x\in\fix(x)$  
is then
$\NC(x) = \{y\in\Fix(f) : D_u(x, \tf) = D_u(y, \tf)\}$.
A fundamental fact is that Nielsen classes of fixed points are
open in $\Fix(f)$ and compact. This allows one to define the 
fixed point index of a Nielsen Class. It also implies that
for a given $f$ the number of (non-empty) Nielsen classes is finite.
A Nielsen class is called \textit{essential} if it has nonzero fixed
 point index.

Central to the theory of persistence of fixed points is the
correspondence of Nielsen classes of homotopic maps.
\begin{theorem}  Assume
$f_0$ is homotopic to $f_1$ via a homotopy $f_t$ and
 $x_0\in\Fix(f_0)$ and $x_1\in\Fix(f_1)$. 
The following are equivalent
\begin{itemize}\label{three2}
\item $D_u(x, \tf) = D_u(y, \tg)$ where $\tf$ and $\tg$ are
equivariantly homotopic.
\item There exists an arc $\gamma$ connecting $x$ to $y$ with
$f_t\circ\gamma $ homotopic to $\gamma$ rel endpoints.
\item There exist lifts $\tx, \ty, \tf, \tg$  to $\tMu$ with 
$\tf$ and $\tg$ equivariantly homotopic and 
$\tf(\tx) = \tx$ and $\tg(\ty) = \ty$.
\end{itemize}
\end{theorem}
If any of these conditions hold we say that $x$ and $y$ \textit{correspond
under the homotopy}. Note that if $x'$ and $y'$ are in the
same Nielsen classes as $x$ and $y$, respectively, then $x'$ and
$y'$ also correspond under the homotopy. Thus it makes sense
to say that $\NC(x,f)$ and $\NC(y,g)$ correspond under the homotopy.
The fundamental theorem of Nielsen theory is
\begin{theorem}
Assume that $x\in\fix(f)$ has an essential Nielsen class. 
If $g$ is homotopic to $f$ there exists a $y\in\fix(g)$ with
$y$ connected to $x$ by the homotopy. Further, the fixed 
point indices of $\NC(x,f)$ and $\NC(y,g)$ are equal.
\end{theorem}

\subsection{periodic points} We say that $x$ is a period
$k$ point for $f$ if $f^k(x) = x$ and this is the
least $k>0$ with this property. This is obviously stronger
than just $x\in\fix(f^k)$.  The complication for persistence
of periodic points under homotopy is that a period $k$ point could have
period $m$ points with $m<k$ which are Nielsen equivalent under $f^k$.
This is an issue because  we
would like the period $k$ point to persist as a period $k$ point.
In the language of
bifurcation theory we need hypothesis on the Nielsen
class of a period point
which disallows its members from undergoing a period dividing bifurcation
under a homotopy. 

Let us first define  displacements of period points.  Fix a lift  $\tfk$
  of $f^k$ and assume $x\in\fix(f^k)$.
If $\tfk(\tx) = \delta \tx$ say  
$\Delta_u(\tx, \tfk) = \delta$ and then for $x\in M$,
$$
D_u(x, \tfk) = \{f_\#^k(\gamma) \delta \gamma^{-1}: \gamma\in \pi_1\}.
$$
A period $k$ orbit is collapsible if its displacement looks like the
 displacement under $\tf^k$ of a periodic point of lesser period.
This leads to the following definition.
\begin{definition}
A period $k$ point for $f$ is  called collapsible (or irreducible) if there
are $\ell>1, m\geq1$ with $\ell m = k$ and an element $\sigma\in\pi_1$ with
$$
f_\#^{(l-1)m}(\sigma)\;
f_\#^{(l-2)m}(\sigma) \dots
f_\#^{m}(\sigma)\; 
\sigma \in D_u(x, \tfk) 
$$
If $x$ is not collapsible, then $x$ and its Nielsen class are
called uncollapsible.
\end{definition}
\begin{theorem}\label{perun} Assume that $x$ is period $k$ and $NC(x,f^k)$ is
essential and uncollapsible. If $g$ is homotopic to $f$ there exists a 
period $k$ point $y$  under $g$  with
$y$ connected to $x$ by the homotopy. Further, the fixed 
point indices of $NC(x,f^k)$ and $NC(y,g^k)$ are equal.
\end{theorem}
In light of this theorem, a period $k$ point and it Nielsen
class are called \textit{persistent} if it is essential and
uncollapsible. Thus the displacement of any persistent point will
be present in any homotopic map.

In general, it is hard to check from the algebra if a given periodic
point is uncollapsible. There is a  sufficient condition
which combined  with
Theorem~\ref{perun}
 gives a simple condition for persistence of period $k$ orbits.
\begin{lemma}\label{different}
If for some $k$, every point in $\fix(f^k)$ is in 
a different Nielsen class under $f^k$, then every period
$k$ point of $f$ is uncollapsible. Thus if
every lift of $f^k$ to $\tMu$ has at most one fixed point, 
then every period $k$ point with nonzero index is
persistent.
\end{lemma}
\begin{proof}
We first prove the contrapositive of the first sentence. Thus
for some $x$, $\ell, m\geq 1$ with $m\ell = k$, we have 
\begin{equation} \label{one} \tf^{\ell m}(\tx) = 
f_\#^{(l-1)m}(\sigma)\; 
f_\#^{(l-2)m}(\sigma) \dots
f_\#^{m} (\sigma)\;
\sigma \tx.
\end{equation}
This implies
\begin{equation}  \tf^{\ell m}(\tf^m(\tx)) = \tf^m(\tf^{\ell m}(\tx)) = 
f_\#^{\ell m}(\sigma)\; 
f_\#^{(l-1)m}(\sigma) \dots
f_\#^{m}(\sigma)\;
\tf^m(\tx).
\end{equation}
On the other hand, \eqref{one} implies 
\begin{equation}
\tf^{\ell m}(\sigma\tx) = f_\#^{\ell m}(\sigma)\tf^{\ell m}\tx = 
f_\#^{\ell m}(\sigma)\; 
f_\#^{(l-1)m}(\sigma) \dots
f_\#^{m}(\sigma)\;
(\sigma \tx).
\end{equation}
Thus the same lift of $f^k$ fixes both $\tf^m(\tx)$ and $\sigma \tx$.
There are two alternatives. If $\tf^m(\tx) = \sigma \tx$,
then $x$ has period $\leq m < k$, a contradiction to the
assumption that $x$ has least period $k$. The other alternative is
that some lift of $f^k$ fixes two points, proving the contrapositive.

For the second sentence, since every lift can fix at most
one point, by Theorem\ref{three2}, no two period $k$ points are Nielsen equivalent
under $f^k$ and so each periodic Nielsen class has at most 
one point. If that point has nonzero index, by Theorem~\ref{perun} 
it is persistent,
\end{proof}

The expansivity condition used in the universal cover
 introduced in Definition~\ref{coverexpand}
allows us to invoke Lemma~\ref{different}.
\begin{theorem}\label{per3} If $f$ is $\tMu$-expanding then for
all $k$, all the points in $\fix(f^k)$ are in different
Nielsen classes. Thus if they have non-zero index, they
are persistent.
\end{theorem}

\subsection{examples}\label{uexpand}
To examine the examples $\phi_1, \phi_2$
and $\phi_3$ note that they all satisfy the following property:
there is a integer $M>1$ so that for any nontrivial word
$w$, $\ell(\psi(w)) = M \ell(w)$ where $\ell$ is the reduced
word length. It follows then that their tight maps made on the wedge 
of circles are $\tMu$-expanding. When $M$ is a compact, connected Riemannian
manifold then any  $f:M\raw M$ that is $\tMu$-expanding
is expanding under Gromov's definition and
thus is topologically conjugate to a infra-nil-endomorphism \cite{gromov2}.  

\subsection{Context}
Lemma~\ref{different} is inspired by a lemma 
from Asimov-Franks  \cite{AF}. For excellent
expositions of  Nielsen theory see the works
of Jiang, for example, \cite{jiang1, jiang2}.

\section{global persistence}\label{GP}
It turns out that in many homotopy classes there is a 
``minimal model'' in the sense that every map in its
class has at least its dynamics. 
\begin{definition}
If $f$ is $\tMu$-expanding and in addition all periodic
orbits have nonzero index and  the set of its periodic
points is dense, then $f$ is called a minimal model in its homotopy
class.
\end{definition}
The terminology minimal model is justified by the following theorem:
\begin{theorem}[Handel]\label{minthm}
Assume $\phi:M\raw M$ is a minimal model and $g$ is homotopic
to $\phi$. Then there exists a compact, $g$-invariant set $Y$ and
a continuous, onto map $\omega:Y\raw M$ which is homotopic to the 
inclusion so that  
\begin{equation*}
 \begin{CD}
 Y @>{g_{\vert Y}}>>Y\\
 @V{\omega}VV    @VV{\omega}V\\
 M @>{\phi}>>M%\\
 \end{CD}
 \end{equation*}
\end{theorem}

\begin{proof}[Sketch] Fix equivariantly
homotopic lifts $\tg$ and $\tphi$ to the
universal cover. Recall that  two orbits $(\tx, \tphi)$ and $(\ty,\tg)$
\textit{$\tMu$-shadow} if there is some
$K$ with  $\td(\tphi^n(\tx), \tg^n(\ty)) < K$ for all $n$ where
$\td$ is some equivariant metric on $\tMu$ which induces the
manifold topology. 

Let $$\tY = \{\ty\in\tMu: \ \text{there exists a}\ (\tx,\tphi) \ \text{that}\ 
\tMu\text{-shadows} (\ty,\tg)\}$$
 and define $\tomega:\tY\raw \tM_u$ by sending
$\ty$ to the $\tx$ it $\tMu$-shadows. Note that because of the
expansion no two points
can $\tMu$-shadow under $\tphi$, so this point $\tx$ is unique and so
$\tomega$ is well-defined. Clearly $\tomega\circ\tg_{\vert\tY} = 
\tphi\circ \tomega$ and $\tomega\circ\delta = \delta\circ\tomega$
for each deck transformation $\delta$
and so $\tomega$ descends to semi-conjugacy $\omega$
on $Y$, the projection of $\tY$.

The crucial uniformity comes from the analog of Theorem~\ref{gs1},
 namely, that 
$(\tx, \tphi)$ and $(\ty,\tg)$ $\tMu$-shadow with some constant
$K$ if and only if 
they $\tMu$-shadow with constant $c/(\lambda-1)$ where
$c=\max\{\td(\tphi(\tx), \tg(\tx)): \tx\in\tM_u\}$. Note that
$c$ is finite since $\tg$ is equivariantly homotopic to $\tphi$.
With a little work this uniform bound implies that $\tY$ is
closed and $\tomega$ is continuous.

The final step uses the persistence of period points. For every
periodic point $x$ of $\phi$ by Theorem~\ref{per3}, there is a periodic
point $y$ of $g$ that is Nielsen equivalent to it. As in Lemma~\ref{pergs}
this implies that they have lifts that $\tMu$-shadow.
Thus the image of the compact set $Y$ contains all the periodic
point of $\phi$. But this set is dense and so $\omega(Y) = M$.
\end{proof}

\begin{remark}\label{tMuhyper}
A homeomorphism $h$ is called
$\tMu$-hyperbolic 
if there a pair of an equivariant pseudo-metrics $\td_1, \td_2$
with $\td = \sqrt{\td_1^2  + \td_2^2}$ a metric inducing the manifold
topology  
and a constant $\lambda>1$ so that $\td_1(\th(\tx), \th(\ty)) 
\geq \lambda \td_1(\tx, \ty)$ and $\td_2(\th\I(\tx), \th\I(\ty)) 
\geq \lambda \td_2(\tx, \ty)$. A homeomorphism is called  a minimal
model if it is $\tMu$-hyperbolic with dense periodic points all
essential.  Theorem~\ref{minthm} is also true in this case.
\end{remark}

\subsection{examples}\label{minmod}
If $A$ is expanding then $\Phi_A:\T^b\raw\T^b$  is a minimal  model. Thus
if $f:\T^b\raw\T^b$ is homotopic to $\Phi_A$, or equivalently
$f_* = A$, then Theorem~\ref{minthm} applies yielding
a semiconjugacy from $f$ acting on a subset of $\T^b$ onto
$\Phi_A$ acting on $\T^b$. On the other hand, $f$ is $H_1$-expanding
and so Theorem~\ref{FS1} applies yielding a semiconjugacy from
$f$ acting on all of $\T^b$ into $\T^b$. We have seen that both
semiconjugacies are equivalent to shadowing and for $M = \T^b$,
$\tMA = \tM_u = \R^b$ we see that the two semiconjugacies
are the same and for the torus the semiconjugacy $\beta$ is
$\T^b$ onto $\T^b$.

As noted in Example~\ref{uexpand}, the maps
 $\phi_1, \phi_2, \phi_3$ are all
$\tMu$ expanding.
Since they each have dense essential periodic
points they are minimal models in their homotopy
classes and so Theorem~\ref{minthm} applies to them all. 
 One useful consequence holds for the map
$\phi_1$ which acts as the identity on homology and thus had
a rotation set as shown in Figure~\ref{rotset}. If $g$ is homotopic to $\phi_1$
then $\rho(\phi_1) \subset \rho(g)$. This follows since $\tMu$-shadowing
implies $\tMa$-shadowing and thus the orbits have the same rotation
vector.

\subsection{mapping torus and suspension flow}\label{Nsusp}
In contrast to homology,
$\pi_1(T_f)$ is ``big enough'' to distinguish different
Nielsen classes of not just fixed points but periodic
points as well. There is one point which requires attention:
two points on same period $k$ orbit
can certainly be in different Nielsen classes under
$f^k$, but their corresponding orbits in the suspension are the
same and thus obviously freely homotopic. So the proper condition
involves the orbits of  the Nielsen class. 
\begin{theorem}[Jiang]\label{jiang3} 
Two period $k$ points satisfy
$\NC(x, f^k)
= \NC(f^i(y), f^k)$ for some $i = 0, \dots, k-1$
if and only if the corresponding closed curves are freely
homotopic in $T_f$.
\end{theorem}

\subsection{example}\label{pione}
Returning to the example of $d:S^1\raw S^1$ via $d(\theta) = 2\theta$,
recall from section~\ref{BFsect} that $H_1(T_d)$ is just $\Z$ and records
just periods of periodic points. 
In contrast, the mapping torus
satisfies $\pi_1(T_d) = \langle z, g: z^{-1} g z = g^2\rangle$. 
This is one of the simpler of the well-studied Baumslag-Solitar groups, 
namely $BS(1,2)$ \cite{BS}. It has a normal subgroup $N$
isomorphic to the group of $2$-adic integers with $BS(1,2)/N \sim \Z$.
It would be interesting to connect the twisted conjugacy classes
of periodic points of $d$ with the algebraic structure of $BS(1,2)$.

\begin{remark} In \cite{Bdamster} the condition in 
Jiang's theorem is incorrectly
stated as  $\NC(x, f^k) = \NC(y, f^k)$ rather than  $\NC(x, f^k)
= \NC(f^i(y), f^k)$. The statement about the relation of strong
Nielsen equivalence to isotopy in the suspension flow needs to
be similarly amended. The statement there about Abelian Nielsen equivalence
and being homologous in the suspension is completely wrong
as the example in Section~\ref{BFsect} shows.
\end{remark}

\subsection{context} Handel's Theorem~\ref{minthm} was originally proved
for \pA\ homeomorphisms.  The proof given above
for Theorem~\ref{minthm} is a sketch of Handel's proof in \cite{handel}.
Jiang's Theorem~\ref{jiang3} is proved in \cite{jiang3}. The 
$\tMu$-displacements of periodic points in the examples on the
wedge of two circles can be computed using the methods of 
Fadell and Husseini \cite{FH}, see also \cite{jiang2}.

\section{Appendix}\label{appendix}
The numerical algorithm for producing the figures of
the image of $\tbeta$ for the examples is based
on two facts. The first is that 
$\tphi(\vzero) = \vzero$, implies using~\eqref{desired} that
 $\tphi^k(\vn) = A^k\vn$ for all
$\vn\in\Z^2$. Thus if $\tz\in\tX$ is such that $\tphi^j(\tz) = \vn$,
then $\tphi^{j+i} (\tz) = A^i \vn$ and so $A^{-(j+i)}\tphi^{j+i}(\tz)
= A^{-j}\vn$ for all $i$. Thus $\tbeta(\tz) = A^{-j}\vn$.
This tells us that an approximation $A^{-k}\tphi^k$ for
finite $k$ will be an accurate computation of $\tbeta(\tz)$
for any $\tz$ with $\tphi^j(\tz) \in\Z^2$ for a $j\leq k$. 

The second fact concerns how to compute the $\tz$ that eventually land
on the lattice under $\tphi$. 
Let $I_a = [0,1]\times \{0\}$ and $I_b = \{0\}\times [0,1]$
and so $W_0=I_a\cup I_b$ is a fundamental domain of $\tX_a$ (with two
extra points). Thus $W_k := A^{-k}\tphi^k(I_a)\cup A^{-k}\tphi^k(I_b)$
is an approximation of $\tbeta(I_a \cup I_b)$. An approximation
of $\tbeta(\tX)$ is obtained by hitting $W_k$ with all the deck transformations
from $\Z^k$. By the previous paragraph, $W_k$ will be accurate, i.e.\
be in $\beta(\tx)$, at those points $\tz\in W_0$ for which $\tphi^j(\tz)\in
\Z^2$ for $ j\leq k$.

To be specific we restrict now to Example~\ref{eg2} and drop subscripts
and so $\psi = \psi_2$ and $\phi = \phi_2$.
Since $\tphi$ is linear with an expansion of $4$ on
pieces of $I_a$ and $I_b$, it is clear that the points 
$\tz\in W_0$ for which $\tphi^j(\tz)\in
\Z^2$ for $ j\leq k$ are 
exactly 
$$
\tz_i = (\frac{i}{4^k}, 0) \in I_a \ \text{and}\ 
\tz_i' = (0, \frac{i}{4^k})\in I_b,
$$
for $0\leq i \leq 4^k$. The value of $\tphi^k$ at these points
can be computed from the words $\psi^k(a)$ and $\psi^k(b)$. Specifically,
in $\R^2$ letting $\kappa(a) = A^{-k}((1,0)^T)$ and  
$\kappa(b) = A^{-k}((0,1)^T)$ we have
$$
\tbeta(\tz_i) = A^{-k} \tphi^k(\tz_i) = 
\sum_{\ell=1}^i \kappa(\psi^k(a)_\ell)\ \text{and} \  
\tbeta(\tz_i') = A^{-k} \tphi^k(\tz_i') = 
\sum_{\ell=1}^i \kappa(\psi^k(b)_\ell) 
$$
Now $\psi^2(a) = a a a b a a \dots$  and
so $\phi^2(5/4^2, 0) = (1,0) + (1,0) + (1,0)+ (0,1) + (1,0) = (3,1)$.
Thus $$
A^{-2}\phi^2(5/4^2,0) = \kappa(a) + \kappa(a) + \kappa(a) + \kappa(b) + 
\kappa(a) = \sum_{\ell=1}^5 \kappa(\psi^2(a)_\ell) = (3/8, -1/8).$$ 
Thus using the result of the first paragraph, $\tbeta(5/4^2, 0) = (3/8, -1/8)$.

Once the images under $\tbeta$ of the points $\tz_i$ and $\tz_i'$
have been computed, they are plotted with line segments drawn
between consecutive images. The resulting image is then accurate 
at ``corners'' but only an approximation elsewhere.

There is a conspicuous inaccuracy visible in Figure~\ref{two_images} (left)
for
Example~\ref{eg2}. We know that $(1/3,0)$ and $(0,1/3)$ are fixed points
of $\phi_2$ that have the same displacements and so they globally
shadow and thus 
have the same $\tbeta$ value. This overlap or intersection
of the two branches does not show up in 
Figure~\ref{two_images} (left). 
The figure is drawn with $k=4$ and the approaching
corners are thus the images of the points $\tz_i$ and $\tz_i'$ which
approximate $(1/3, 0)$ and $(0,1/3)$. It is clear that in
this context, $1/3$ not well approximated by any $\frac{i}{4^k}$.

The intersections visible in Figure~\ref{two_images} (left)
 come from another source, namely,
the non-injectivity of $\tphi_2$ acting on $\tX_a$. For if 
$\tphi_2^k(\tx) = \tphi_2^k(\ty)$ certainly for all $m>0$
$A^{-m-k}\tphi_2^{m+k}(\tx) = A^{-m-k}\tphi_2^{m+k}(\ty)$
and so $\tbeta(\tx) = \tbeta(\ty)$. For example, 
$\tphi_2(1/2,0) = \tphi_2(1, -1/2) = (2,0)$ by direct calculation.

\bibliographystyle{plain}
\bibliography{simons}

\end{document}